\documentclass{amsart}
\usepackage{amssymb,latexsym,subfigure,curves,epic}

\numberwithin{equation}{section}

\newtheorem{theorem}{Theorem}[section]
\newtheorem{proposition}[theorem]{Proposition}
\newtheorem{corollary}[theorem]{Corollary}

\newtheorem{lemma}[theorem]{Lemma}

\newtheorem{definition}[theorem]{Definition}
\newtheorem{algorithm}[theorem]{Algorithm}

\def\proof{\smallskip\noindent {\bf Proof. }}
\def\endproof{\hfill$\square$\medskip}
\def\dc{\lambda}
\def\da{\tilde{\lambda}}
\def\mc{\mathcal{C}}
\def\ma{\mathcal{A}}

\newcommand{\journal}[1]{{\sl #1}}
\newcommand{\arttitle}[1]{{\rm #1}}

\newcommand{\subseqnoref}[1]{\medskip \noindent {\bf #1.}}

\begin{document}{}

\title[Finite posets and Ferrers shapes]{{\ }\\[-.9in]
Finite posets and Ferrers shapes}

\author{Thomas Britz}
\address{\noindent Department of Mathematics, University of 
  \AA rhus, 8000 \AA rhus, Denmark}
\email{britz@imf.au.dk}

\author{Sergey Fomin}
\address{Department of Mathematics, Massachusetts Institute of
  Technology, Cambridge, MA 02139, USA.
\emph{Current address:} Department of Mathematics,  
University of Michigan, 
Ann Arbor, MI 48109, USA}
\email{fomin@math.mit.edu, fomin@math.lsa.umich.edu}
  
\thanks{The second author was supported in part 
by NSF grant \#DMS-9700927.}

\date{December 14, 1999}

\maketitle

\tableofcontents

\vspace{-.5in}

\begin{section}{Introduction}
\label{sec:intro}

This survey, written at the suggestion of G.-C.~Rota, 
focuses on the fundamental correspondence---originally discovered 
by C.~Greene~\cite{greene3}, 
following his joint work with D.~J.~Kleitman~\cite{GK1}---that associates 
a Ferrers shape $\dc(P)$ to every finite 
poset~$P$. 
The number of boxes in the first $k$ rows (resp.\ columns) of $\dc(P)$
equals the maximal number of elements in a union of $k$ chains
(resp.\ antichains) in~$P$. 

The correspondence $P\mapsto \dc(P)$ is intimately
related to at least three areas of discrete mathematics:
combinatorial optimization, lattice theory, and the
combinatorics of tableaux.
In this article, we bring together the main results in the subject, 
along with complete proofs. 

The structure of the paper is as follows. 
In Section~\ref{sec:main}, we state the main theorems.  
Sections~\ref{sec:schensted}--\ref{sec:schutzenberger}
are devoted to tableau-theoretic applications.
In Section~\ref{sec:sat}, the results on saturated families of chains
and antichains are derived. 
Section~\ref{sec:nilpotent} discusses an interpretation of the
main correspondence, 
due to E.~R.~Gansner and M.~Saks, 
in terms of sizes of Jordan blocks of nilpotent elements in the incidence
algebra of~$P$. 

Sections~\ref{sec:flow-prelim}--\ref{sec:gansner34} are devoted to proofs. 
We begin by reproducing A.~Frank's remarkable proof~\cite{frank} 
of the main ``duality theorem'' (Theorem~\ref{thm:dual}) that uses
network flows. 
We then provide three proofs of the ``monotonicity
theorem'' (Theorem~\ref{thm:mono}):
a new beautiful lattice-theoretic proof contributed by C.~Greene
(reproduced with permission); 
E.~Gansner's amazingly short proof~\cite{gansner}  
utilizing the nilpotent-matrices interpretation mentioned above; 
and a proof based on Frank's approach,
which as a byproduct yields an augmenting-path result for maximal
chain families (Theorem~\ref{thm:fomin1}). 
The latter proof, as well as our proofs of Theorems
\ref{thm:fomin3} and~\ref{thm:fominG}, are new, although some of the
ingredients were recycled 
from~\cite{fomin1, fomin3, fomin4, frank, GK1}.   

\end{section}

\begin{section}{Main Theorems}
\label{sec:main}

Let $P$ be a finite partially ordered set of cardinality~$n$. 
A \emph{chain} is a totally ordered subset of~$P$. 
An \emph{antichain} is a subset of $P$ in which no two elements are 
comparable. 
The famous theorem of Dilworth~\cite{dilworth1} states that the
maximal size of an antichain in $P$ is equal to the minimal number of
chains into which $P$ can be partitioned. 
This theorem has an easy ``dual'' counterpart, 
in which the words ``chain'' and ``antichain'' are interchanged
(see~\cite{mirsky} or \cite[Theorem~6.2]{van-lint-wilson}). 

Dilworth's theorem and its dual have a beautiful and powerful common
generalization due to Curtis Greene (Theorem~\ref{thm:dual} below).
 
For $k=0,1,2,\dots$, let $a_k$ (resp.~$c_k$) denote the maximal
cardinality of a union of $k$ antichains (resp.\ chains). 
Let $\dc_k=c_k-c_{k-1}$ and $\da_k=a_k-a_{k-1}$ 
for all $k\geq 1$.

\begin{theorem}[The Duality Theorem for Finite Partially Ordered Sets] 
\label{thm:dual} 
\ \\
The sequences 
$\dc=(\dc_1, \dc_2, \dots)$ and 
$\da=(\da_1, \da_2, \dots)$ 
are weakly decreasing, and form conjugate partitions of the
number~$n=|P|$. 
\end{theorem}

This theorem was first obtained by
C.~Greene~\cite{greene3} as a corollary of 
another prominent result, due to C.~Greene and
D.~Kleitman~\cite{GK1}. 
A few years later, it was rediscovered and given an alternative proof
in~\cite{fomin1}. 
Other proofs appeared as well; 
we would like to single out an elegant proof given
by A.~Frank (reproduced in~\cite{engel, felsner} and in 
Section~\ref{sec:frank's-proof} below). 
 
The Duality Theorem associates to every finite poset the Young diagram
(or Ferrers shape) whose row lengths are $\dc_1, \dc_2,\dots$,
and whose column lengths are $\da_1, \da_2, \dots$. 
We will identify this shape with the partition~$\dc$ and denote it 
by~$\dc(P)$. 
 
To illustrate, consider the poset $P$ in Figure~\ref{fig:dual}.
For this poset, we have $c_0=0$, $c_1=4$, $c_2=c_3=\cdots=6$, 
implying $\dc=(4,2)$, 
while $a_0=0$, $a_1=2$, $a_2=4$, $a_3=5$,
$a_4=a_5=\cdots=6$, implying that $\da=(2,2,1,1)$, a shape conjugate to~$\dc$. 

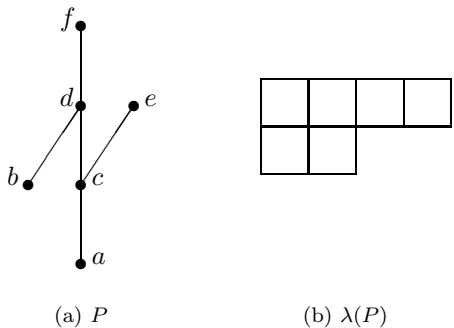
\begin{figure}[ht]
\centering
\setlength{\unitlength}{2pt}
\mbox{
    \subfigure[$P$]{
        \begin{picture}(40,45)(0,10)
        \put(20,10){\line(0,1){45}}
        \put(10,25){\line(2,3){10}}
        \put(20,25){\line(2,3){10}}
        \multiput(20,10)(0,15){4}{\circle*{2}}
        \put(10,25){\circle*{2}}
        \put(30,40){\circle*{2}}
        \put(22,10){\makebox{$a$}}
        \put(6,25){\makebox{$b$}}
        \put(22,25){\makebox{$c$}}
        \put(16,40){\makebox{$d$}}
        \put(32,40){\makebox{$e$}}
        \put(16,55){\makebox{$f$}}
        \end{picture}
        } 
\qquad
    \subfigure[$\dc(P)$]{
        \begin{picture}(30,20)(10,10)
\setlength{\unitlength}{3pt}
        \put(6,18){\line(1,0){12}}
        \multiput(6,24)(0,6){2}{\line(1,0){24}}
        \multiput(6,18)(6,0){3}{\line(0,1){6}}
        \multiput(6,24)(6,0){5}{\line(0,1){6}}
        \end{picture}
        } 
}
\caption{The Duality Theorem}
\label{fig:dual}
\end{figure} 

As an immediate corollary of Theorem~\ref{thm:dual}, 
the number of 
rows in $\dc=\dc(P)$ is equal to $\da_1\,$, 
a reformulation of Dilworth's theorem. 

Various attempts have been made
(see, e.g.,~\cite{west1,linial,felsner,HSH}) 
to generalize Theorem~\ref{thm:dual} to directed graphs. 
In this survey, we do not discuss these generalizations. 


The Duality Theorem naturally associates a Ferrers shape to any finite
poset. 
The following result shows that this correspondence is, in some sense,
``functorial.'' 

\begin{theorem}[The Monotonicity Theorem] 
\label{thm:mono} \cite{fomin1} 
Let $p$ be a maximal (or minimal) element of a finite poset~$P$. 
Then $\dc(P-\{p\})\subset \dc(P)$.
\end{theorem}

For example, the poset $P$ in Figure~\ref{fig:dual} has maximal
elements $e$ and~$f$. 
The shapes $\dc(P-\{e\})$ and $\dc(P-\{f\})$ are shown in 
Figure~\ref{fig:mono}; both are contained in~$\dc(P)$. 

\begin{figure}[ht]
\centering
\setlength{\unitlength}{3pt}
    \subfigure[$\dc(P\!-\!\{e\})$]{
        \begin{picture}(30,15)(-3,0)
\setlength{\unitlength}{3pt}
        \put(0,0){\line(1,0){6}}
        \multiput(0,6)(0,6){2}{\line(1,0){24}}
        \multiput(0,0)(6,0){2}{\line(0,1){6}}
        \multiput(0,6)(6,0){5}{\line(0,1){6}}
        \end{picture}
        } 
\qquad\qquad
    \subfigure[$\dc(P\!-\!\{f\})$]{
        \begin{picture}(30,15)(-5,0)
\setlength{\unitlength}{3pt}
        \put(0,0){\line(1,0){12}}
        \multiput(0,6)(0,6){2}{\line(1,0){18}}
        \multiput(0,0)(6,0){3}{\line(0,1){6}}
        \multiput(0,6)(6,0){4}{\line(0,1){6}}
        \end{picture}
        } 
\caption{The Monotonicity Theorem}
\label{fig:mono}
\end{figure}
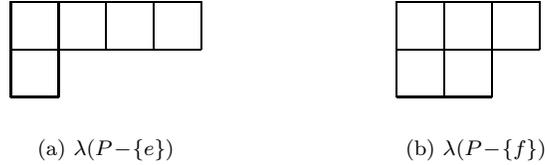 

In the Monotonicity Theorem, the restriction for $p\in P$ to be an
extremal element cannot be dropped. 
A counterexample is given in Figure~\ref{fig:notmono}. 

\begin{figure}[ht]
\centering
\setlength{\unitlength}{3pt}
\mbox{
\qquad
    \subfigure[$P$]{
\setlength{\unitlength}{2pt}
        \begin{picture}(50,40)(-8,10)
        \put(10,25){\line(2,3){10}}
        \put(20,10){\line(0,1){30}}
        \put(20,10){\line(2,3){10}}
        \multiput(10,25)(10,0){3}{\circle*{2}}
        \put(20,10){\circle*{2}}
        \put(20,40){\circle*{2}}
        \put(16,25){\makebox{$p$}}
        \end{picture}
        } 
    \subfigure[$\dc(P)$]{
        \begin{picture}(30,25)(0,12)
\setlength{\unitlength}{3pt}
        \multiput(6,12)(0,6){2}{\line(1,0){6}}
        \multiput(6,24)(0,6){2}{\line(1,0){18}}
        \multiput(6,12)(6,0){2}{\line(0,1){18}}
        \multiput(18,24)(6,0){2}{\line(0,1){6}}
        \end{picture}
        } 
\qquad 
    \subfigure[$\dc(P\!-\!\{p\})$]{
        \begin{picture}(30,25)(-3,15)
\setlength{\unitlength}{3pt}
        \multiput(6,18)(0,6){3}{\line(1,0){12}}
        \multiput(6,18)(6,0){3}{\line(0,1){12}}
        \end{picture}
        } 
}
\caption{A counterexample: $\dc(P-\{p\})\not\subset \dc(P)$}
\label{fig:notmono}
\end{figure}
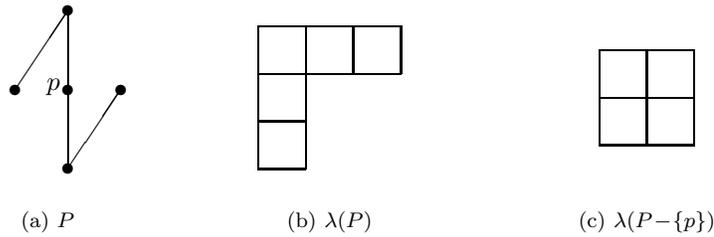

Theorem~\ref{thm:mono} implies that any linear extension 
$\varphi :P\rightarrow [n]=\{1,\dots,n\}$ 
of~$P$ gives rise to a standard Young tableau $T$ of shape
$\dc(P)$ (see ~\cite[p.~312]{stanley2}) 
defined by the condition that the entries 
$1,\dots,k$ of~$T$ form the shape $\dc(\varphi([1,k]))$. 
As an example, consider the poset in Figure~\ref{fig:dual} 
and its linear extension given by $\varphi(a)=1$, $\varphi(b)=2$,
\dots, $\varphi(f)=6$. 
The resulting standard tableau is given in Figure~\ref{fig:tableau}.

\begin{figure}[ht]
\centering
\mbox{
\setlength{\unitlength}{1.9pt}
\begin{picture}(40,55)(0,10)
\put(20,10){\line(0,1){45}}
\put(10,25){\line(2,3){10}}
\put(20,25){\line(2,3){10}}
\multiput(20,10)(0,15){4}{\circle*{2}}
\put(10,25){\circle*{2}}
\put(30,40){\circle*{2}}
\put(22,10){\makebox{$1$}}
\put(6,25){\makebox{$2$}}
\put(22,25){\makebox{$3$}}
\put(16,40){\makebox{$4$}}
\put(32,40){\makebox{$5$}}
\put(16,55){\makebox{$6$}}
\end{picture}

\qquad\qquad

\setlength{\unitlength}{3pt}
\begin{picture}(30,12)(0,-10)
\put(0,0){\line(1,0){12}}
\multiput(0,6)(0,6){2}{\line(1,0){24}}
\multiput(0,0)(6,0){3}{\line(0,1){6}}
\multiput(0,6)(6,0){5}{\line(0,1){6}}
\put(3, 9){\makebox(0,0){$1$}}
\put(9, 9){\makebox(0,0){$3$}}
\put(15,9){\makebox(0,0){$4$}}
\put(21,9){\makebox(0,0){$6$}}
\put(3, 3){\makebox(0,0){$2$}}
\put(9, 3){\makebox(0,0){$5$}}
\end{picture}
}
\caption{A linear extension and the associated standard tableau}
\label{fig:tableau}
\end{figure}
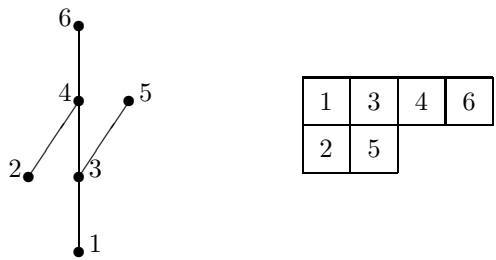

A generalization of Theorem~\ref{thm:mono} to path families in acyclic
directed graphs 
was given by E.~R.~Gans\-ner \cite[Theorem~3.2]{gansner}. 

We will next describe a simple 
recursive algorithm for computing the shapes $\dc(P')$ for all order
ideals~$P'$ of a given finite poset~$P$. 

\begin{theorem}[Recursive computation of the shape]
\label{thm:fomin3} 
\cite{fomin4} 
Let $p_1,\dots,p_k$ be the full list of maximal elements in~$P$. 
Then the shape $\dc=\dc(P)$ is uniquely determined by the shapes 
$\dc(P-\{p_1\})$, \dots, $\dc(P-\{p_k\})$, as follows. 

If $\dc(P-\{p_1\})=\cdots=\dc(P-\{p_k\})=\dc'$, then $\dc$ is  
obtained by adding a box into the $k$'th row of~$\dc'$.
Otherwise, $\dc=\bigcup_i \dc(P-\{p_i\})$, by Theorem~\ref{thm:mono}. 
\end{theorem}

Theorem~\ref{thm:fomin3} can be used to recursively
compute the shapes $\dc(P')$ for all order ideals $P'$ of~$P$; 
such ideals form a distributive lattice denoted~$J(P)$. 
To illustrate, consider the poset $P$ in Figure~\ref{fig:dual}. 
The lattice $J(P)$ is shown in Figure~\ref{fig:lattice}. 
The shapes associated to the elements of~$J(P)$ are computed bottom-up
using the rule of Theorem~\ref{thm:fomin3}.
For example, the element $\{a,b,c\}\in J(P)$ covers 
$\{a,b\}$ and $\{a,c\}$. Since $\dc(\{a,b\})\neq \dc(\{a,c\})$,
we have $\dc(\{a,b,c\})=\dc(\{a,b\})\cup\dc(\{a,c\})$.
On the other hand, the element $\{a,b,c,d,e\}$ covers two elements 
$\{a,b,c,d\}$ and $\{a,b,c,e\}$ 
(obtained by removing maximal elements $e$ and $d$, respectively) 
such that $\dc(\{a,b,c,d\})= \dc(\{a,b,c,e\})=
\setlength{\unitlength}{4.5pt}
\begin{picture}(4,2)(0,0.5)
                        \put(0,0){\line(1,0){1}}
                        \multiput(2,1)(1,0){2}{\line(0,1){1}}
                        \multiput(0,0)(1,0){2}{\line(0,1){2}}
                        \multiput(0,1)(0,1){2}{\line(1,0){3}}
                        \end{picture}$.
Thus the shape $\dc(\{a,b,c,d,e\})$ is obtained by adding a box
into the second row of \hspace{.02in}
$\setlength{\unitlength}{4.5pt}
\begin{picture}(4,2)(0,0.5)
                        \put(0,0){\line(1,0){1}}
                        \multiput(2,1)(1,0){2}{\line(0,1){1}}
                        \multiput(0,0)(1,0){2}{\line(0,1){2}}
                        \multiput(0,1)(0,1){2}{\line(1,0){3}}
                        \end{picture}$.

\begin{figure}[ht]
\centering
\setlength{\unitlength}{4.5pt}
\mbox{
    \subfigure[The order ideals\dots]{
        \begin{picture}(36,37)(0,0) 
        \multiput(15, 5)(-5,5){2}{\line(1, 1){15}}
        \multiput(15,25)(-5,5){2}{\line(1, 1){ 5}}
        \multiput(10,10)( 5,5){2}{\line(1,-1){ 5}}
        \multiput(10,30)( 5,5){2}{\line(1,-1){15}}
        \multiput(15, 5)( 5,5){4}{\circle*{1}}
        \multiput(10,10)( 5,5){4}{\circle*{1}}
        \multiput(15,25)( 5,5){2}{\circle*{1}}
        \multiput(10,30)( 5,5){2}{\circle*{1}}
        \put(15,5.5){\line( 1,1){5}}
        \put(20,10.5){\line(-1,1){5}}
        \put(15,15.5){\line( 1,1){5}}
        \put(20,20.5){\line(-1,1){5}}
        \put(15,25.5){\line( 1,1){5}}
        \put(20,30.5){\line(-1,1){5}}
        \put(14.5, 1){\makebox{$\phi$}}
        \put(7, 10){\makebox{$b$}}
        \put(22,10){\makebox{$a$}}
        \put(11,15){\makebox{$ab$}}
        \put(27,15){\makebox{$ac$}}
        \put(15,20){\makebox{$abc$}}
        \put(32,20){\makebox{$ace$}}
        \put(8, 25){\makebox{$abcd$}}
        \put(26,25){\makebox{$abce$}}
        \put(3, 30){\makebox{$abcdf$}}
        \put(21,30){\makebox{$abcde$}}
        \put(12,37){\makebox{$abcdef$}}
        \end{picture} 
    }
\ 
    \subfigure[\dots and their shapes]{
        \begin{picture}(36,37)(0,0) 
        \multiput(15, 5)(-5,5){2}{\line(1, 1){15}}
        \multiput(15,25)(-5,5){2}{\line(1, 1){ 5}}
        \multiput(10,10)( 5,5){2}{\line(1,-1){ 5}}
        \multiput(10,30)( 5,5){2}{\line(1,-1){15}}
        \multiput(15, 5)( 5,5){4}{\circle*{1}}
        \multiput(10,10)( 5,5){4}{\circle*{1}}
        \multiput(15,25)( 5,5){2}{\circle*{1}}
        \multiput(10,30)( 5,5){2}{\circle*{1}}
        \put(14.5,1){\makebox{$\phi$}}
        \put(6, 10){\makebox{
                        \begin{picture}(1,1)(0,0)
                        \multiput(0,0)(0,1){2}{\line(1,0){1}}
                        \multiput(0,0)(1,0){2}{\line(0,1){1}}
                        \end{picture}}}
        \put(21,10){\makebox{
                        \begin{picture}(1,1)(0,0)
                        \multiput(0,0)(0,1){2}{\line(1,0){1}}
                        \multiput(0,0)(1,0){2}{\line(0,1){1}}
                        \end{picture}}}
        \put(11,14){\makebox{
                        \begin{picture}(1,2)(0,0)
                        \multiput(0,0)(0,1){3}{\line(1,0){1}}
                        \multiput(0,0)(1,0){2}{\line(0,1){2}}
                        \end{picture}}}
        \put(26,14){\makebox{
                        \begin{picture}(2,2)(0,0)
                        \multiput(0,1)(0,1){2}{\line(1,0){2}}
                        \multiput(0,1)(1,0){3}{\line(0,1){1}}
                        \end{picture}}}
        \put(14,20){\makebox{
                        \begin{picture}(2,2)(0,0)
                        \put(0,0){\line(1,0){1}}
                        \put(2,1){\line(0,1){1}}
                        \multiput(0,0)(1,0){2}{\line(0,1){2}}
                        \multiput(0,1)(0,1){2}{\line(1,0){2}}
                        \end{picture}}}
        \put(30.5,20){\makebox{
                        \begin{picture}(3,2)(0,0)
                        \multiput(0,1)(1,0){4}{\line(0,1){1}}
                        \multiput(0,1)(0,1){2}{\line(1,0){3}}
                        \end{picture}}}
        \put(8, 25){\makebox{
                        \begin{picture}(3,2)(0,0)
                        \put(0,0){\line(1,0){1}}
                        \multiput(2,1)(1,0){2}{\line(0,1){1}}
                        \multiput(0,0)(1,0){2}{\line(0,1){2}}
                        \multiput(0,1)(0,1){2}{\line(1,0){3}}
                        \end{picture}}}
        \put(26,25){\makebox{
                        \begin{picture}(3,2)(0,0)
                        \put(0,0){\line(1,0){1}}
                        \multiput(2,1)(1,0){2}{\line(0,1){1}}
                        \multiput(0,0)(1,0){2}{\line(0,1){2}}
                        \multiput(0,1)(0,1){2}{\line(1,0){3}}
                        \end{picture}}}
        \put(3, 30){\makebox{
                        \begin{picture}(4,2)(0,0)
                        \put(0,0){\line(1,0){1}}
                        \multiput(2,1)(1,0){3}{\line(0,1){1}}
                        \multiput(0,0)(1,0){2}{\line(0,1){2}}
                        \multiput(0,1)(0,1){2}{\line(1,0){4}}
                        \end{picture}}}
        \put(21,30){\makebox{
                        \begin{picture}(3,2)(0,0)
                        \put(0,0){\line(1,0){2}}
                        \put(3,1){\line(0,1){1}}
                        \multiput(0,0)(1,0){3}{\line(0,1){2}}
                        \multiput(0,1)(0,1){2}{\line(1,0){3}}
                        \end{picture}}}
        \put(12,37){\makebox{
                        \begin{picture}(4,2)(0,0)
                        \put(0,0){\line(1,0){2}}
                        \multiput(3,1)(1,0){2}{\line(0,1){1}}
                        \multiput(0,0)(1,0){3}{\line(0,1){2}}
                        \multiput(0,1)(0,1){2}{\line(1,0){4}}
                        \end{picture}}}
        \end{picture}
    }
} 
\caption{Recursive computation along the lattice~$J(P)$} 
\label{fig:lattice} 
\end{figure}
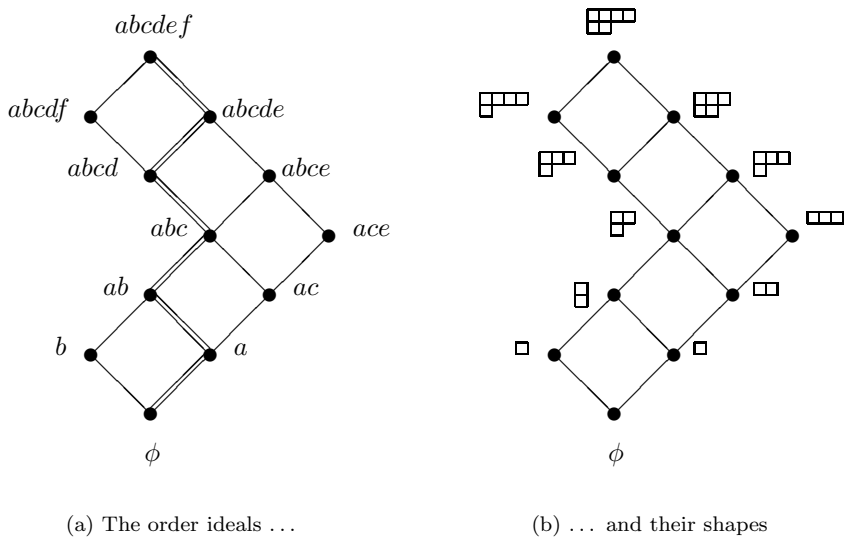



Note that this algorithm can be used as an alternative recursive
definition of the correspondence $P\mapsto \dc(P)$. 
However, if one uses this definition, without invoking the Duality
Theorem, then a natural question arises: why does this recursive procedure
never break down? 
To rephrase, why does each recursive step produce a legal shape of the
right number of boxes? 
There must be a way to answer these questions directly;
this might lead to yet another independent proof of the Duality
Theorem. 


The following theorem, which will prove to be useful in
tableau-theoretic applications of Section~\ref{sec:schensted}, 
provides a more detailed information regarding the growth of the shape
$\dc(P)$ as we add/remove extremal elements to/from~$P$.   


\begin{theorem} 
\label{thm:fominG} 
\cite{fomin3}
Assume that $p_1$ and $p_2$ are extremal
(i.e., maximal or minimal) elements of~$P$,
and suppose that 
$\dc (P-\{ p_1\})=\dc (P-\{ p_2\})$.
Denote $\dc=\dc(P)$, and let the boxes $A$ and $B$ be defined by
\begin{eqnarray}
\label{eq:conditionsG} 
\begin{array}{l}
\dc(P-\{p_1\})=\dc (P-\{ p_2\})=\dc-\{B\} \,, \\[.1in]
\dc(P-\{p_1,p_2\}) =\dc-\{A,B\}\,.
\end{array}
\end{eqnarray}
If $p_1$ and $p_2$ are both maximal or both minimal, 
then $A$ is located either in the same column as $B$ or to the right
of~$B$. 
If, on the other hand, $p_1$ is maximal while $p_2$ is minimal (or
vice versa), then $A$ is either in the same column as $B$ or in the
column immediately to the left of~$B$.  
(See Figure~\ref{fig:fominG}.) 
\end{theorem}


\begin{figure}[ht]
\centering
\setlength{\unitlength}{2.5pt}
\mbox{
    \subfigure[$p_1$ and $p_2$ both minimal (or both maximal)]{  
        \begin{picture}(48,36)(-3,0)
        \put(0, 0){\line(0,1){36}}
        \put(0,36){\line(1,0){42}}
        \put(24,18){\line(1,0){18}}
        \put(24,12){\line(0,1){12}}
        \put(30,12){\line(0,1){24}}
        \multiput(24,12)(0,12){2}{\line(1,0){6}}
        \put(27,15){\makebox(0,0){$B$}} 
        \multiput(24,18)(0,1.5){4}{\line(1,0){18}} 
        \multiput(30,24)(0,1.5){8}{\line(1,0){12}} 
        \end{picture}
    }
    \subfigure[$p_1$ minimal, $p_2$ maximal (or vice versa)]{\qquad  
        \begin{picture}(72,36)(-12,0)
        \put(0, 0){\line(0,1){36}}
        \put(0,36){\line(1,0){42}}
        \put(18,18){\line(1,0){12}}
        \put(24,0){\line(0,1){24}}
        \put(30,12){\line(0,1){12}}
        \multiput(24,12)(0,12){2}{\line(1,0){6}}
        \put(27,15){\makebox(0,0){$B$}} 
        \multiput(18, 0)(1.5,0){4}{\line(0,1){18}} 
        \multiput(24,18)(1.5,0){4}{\line(0,1){6}} 
        \end{picture}
    }
} 
\caption{Theorem~\ref{thm:fominG}: allowable locations of~$A$,
  given~$B$} 
\label{fig:fominG} 
\end{figure}
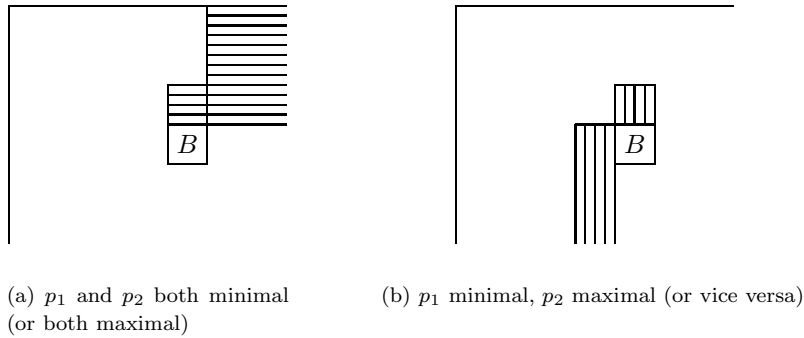

Various subcases of Theorem~\ref{thm:fominG} are exemplified in
Figure~\ref{fig:Gexample}. 
Deleting each of the extremal elements 
$p_1$, $p_2$, $p'_2$, $p''_2$, $\tilde p_2$ from $P$ 
results in the removal of the box~$B$ from $\dc=\dc(P)$. 
Furthermore, deleting $p_1$ together with 
$p_2$ (resp.\ $p'_2$, $p''_2$, $\tilde p_2$) 
results in removing $B$ together with~$A$
(resp.\ $A'$, $A''$, $\tilde A$),
in agreement with Theorem~\ref{thm:fominG}. 

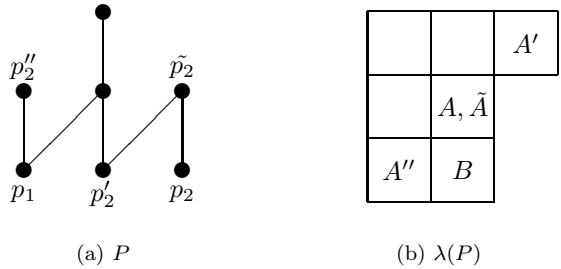
\begin{figure}[ht]
\centering
\setlength{\unitlength}{3pt}
\mbox{
    \subfigure[$P$]{
        \centering
        \begin{picture}(20,27)(0,6) 
        \multiput(0,10)(20,0){2}{\line(0,1){10}} 
        \multiput(0,10)(10,0){2}{\line(1,1){10}} 
        \put(10,10){\line(0,1){20}} 
        \multiput(0,10)(10,0){3}{\circle*{2}} 
        \multiput(0,20)(10,0){3}{\circle*{2}} 
        \put(10,30){\circle*{2}} 
        \put(0 ,7){\makebox(0,0){$p_1$}} 
        \put(10,7){\makebox(0,0){$p_2 '$}} 
        \put(20,7){\makebox(0,0){$p_2$}} 
        \put(0 ,23){\makebox(0,0){$p_2 ''$}} 
        \put(20,23){\makebox(0,0){$\tilde{p_2}$}} 
        \end{picture} 
    }
\qquad \qquad  \qquad 
    \subfigure[$\dc(P)$]{
        \centering
        \begin{picture}(18,18)(0,0) 
\setlength{\unitlength}{4pt}
        \multiput(0,  0)(0,6){2}{\line(1,0){12}} 
        \multiput(0, 12)(0,6){2}{\line(1,0){18}} 
        \multiput(0,  0)(6,0){3}{\line(0,1){18}} 
        \put(18,12){\line(0,1){6}} 
        \put(3,  3){\makebox(0,0){$A''$}} 
        \put(9,  3){\makebox(0,0){$B$}} 
        \put(9,  9){\makebox(0,0){$A, \tilde{A}$}}
        \put(15,15){\makebox(0,0){$A'$}} 
        \end{picture}
    }
} 
\caption{An example illustrating Theorem~\ref{thm:fominG}} 
\label{fig:Gexample} 
\end{figure}

Theorem~\ref{thm:fominG} is sharp~\cite{fomin3} in the sense that for
any shape $\dc$ and any boxes $A$ and~$B$ located in compliance 
with the rules of Figure~\ref{fig:fominG}
(also, $\dc-\{B\}$ and $\dc-\{A,B\}$ should be valid shapes of
$|\dc|-1$ and $|\dc|-2$ boxes, respectively), 
one can produce a poset $P$ together with extremal
elements $p_1$ and $p_2$ of appropriate kind
so that (\ref{eq:conditionsG}) are satisfied. 
The proof of this simple assertion is given at the end of 
Section~\ref{sec:proof-fominG}. 

Our list of known general restrictions governing the growth of the
shape $\dc(P)$ would not be complete without the following simple
result due to E.~R.~Gansner \cite[Theorem~3.4]{gansner}.  
(Gansner's theorem actually holds in greater generality, 
for arbitrary acyclic directed graphs.)

\begin{theorem}
\label{thm:gansner34}
\cite{gansner}
Assume that $p_1$ is a maximal element in $P$,
while $p_2$ is a maximal element in $P-\{ p_1\}$
such that $p_1$ covers $p_2\,$. 
Let the boxes $A$ and $B$ be defined by
\begin{eqnarray}
\label{eq:conditions34}
\begin{array}{l}
\dc(P-\{p_1\})=\dc(P)-\{B\} \,, \\[.1in]
\dc(P-\{p_1,p_2\}) =\dc-\{A,B\}\,.
\end{array}
\end{eqnarray}
Then $A$ is located to the left of $B$. 
(See Figure~\ref{fig:gansner34}.)
\end{theorem}

\begin{figure}[ht]
\centering
\setlength{\unitlength}{2.5pt}
\begin{picture}(48,36)(-3,0)
\put(0, 0){\line(0,1){36}}  
\put(0,36){\line(1,0){42}}  
\put(18,18){\line(1,0){12}} 
\put(0,12){\line(1,0){30}} 
\put(30,12){\line(0,1){6}}
\put(27,15){\makebox(0,0){$B$}}
\multiput(18,0)(1.5,0){5}{\line(0,1){18}}
\multiput(0,0)(1.5,0){12}{\line(0,1){12}}
\end{picture}
\caption{Theorem~\ref{thm:gansner34}: allowable locations of~$A$,
  given~$B$}
\label{fig:gansner34}
\end{figure}
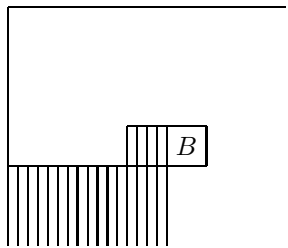
        
To illustrate, consider Figure~\ref{fig:34example}. 
Deleting the maximal element $p_1$ from $P$ 
results in the removal of the box~$B$. 
Subsequent removal of $p_2$ (resp. $p_2 '$) 
results in the removal of the box $A$ (resp. $A'$), 
in compliance with Theorem~\ref{thm:gansner34}. 

\begin{figure}[ht]
\centering
\setlength{\unitlength}{3pt}
\mbox{
    \subfigure[$P$]{
        \centering
        \begin{picture}(20,27)(0,6)
        \put(10,10){\line(0,1){20}}
        \put(20,10){\line(0,1){10}}  
        \multiput(0 ,20)(10,-10){2}{\line(1,1){10}}  
        \multiput(10,10)(10,0){2}{\circle*{2}}
        \multiput(0 ,20)(10,0){3}{\circle*{2}}
        \put(10,30){\circle*{2}}
        \put( 7,31){\makebox(0,0){$p_1$}}
        \put(-3,21){\makebox(0,0){$p_2$}}
        \put(13,21){\makebox(0,0){$p_2 '$}}
        \end{picture}
    }
\qquad \qquad  \qquad
    \subfigure[$\dc(P)$]{
        \centering
        \begin{picture}(18,18)(0,0)
\setlength{\unitlength}{4pt}
        \put(0,  0){\line(1,0){6}}
        \put(0,  6){\line(1,0){12}} 
        \put(12, 6){\line(0,1){12}} 
        \put(18,12){\line(0,1){6}} 
        \multiput(0, 12)(0,6){2}{\line(1,0){18}}
        \multiput(0,  0)(6,0){2}{\line(0,1){18}}
        \put(3,  3){\makebox(0,0){$A$}}
        \put(9,  9){\makebox(0,0){$A'$}}
        \put(15,15){\makebox(0,0){$B$}}
        \end{picture}
    }
}
\caption{An example illustrating Theorem~\ref{thm:gansner34}} 
\label{fig:34example}
\end{figure}
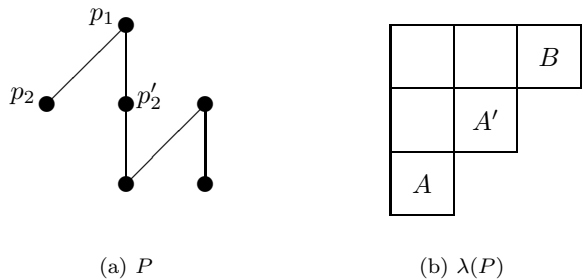


\end{section}


\begin{section}{The Robinson-Schensted correspondence} 
\label{sec:schensted} 

In this section, we explain how the theory of the Robinson-Schensted 
correspondence can be constructed on the poset-theoretic basis laid
out in Section~\ref{sec:main}. 
Our presentation closely follows that of~\cite{fomin3} 
(see also~\cite{fomin6} and~\cite{leeuwen1,roby}), the main tool
being Theorem~\ref{thm:fominG}. 


Let $\sigma =(\sigma(1),\dots,\sigma(n))$ be a permutation
of~$[n]=\{1,\dots,n\}$. 
The associated \emph{permutation poset} $P_\sigma$ is 
the set of ordered
pairs $(i,\sigma (i))$, $i=1,\dots,n$, with the partial order induced
from the product of chains $[n]\times [n]$: 
\[
(i,\sigma(i))\leq(j,\sigma(j))\,\Longleftrightarrow
\text{$i\leq j$ and $\sigma(i)\leq \sigma(j)$} \,. 
\]
An example is given in Figure~\ref{fig:permposet}. 
Note that the poset $P_{412563}$ in Figure~\ref{fig:permposet}b is
isomorphic to the poset~$P$ in Figure~\ref{fig:dual}. 

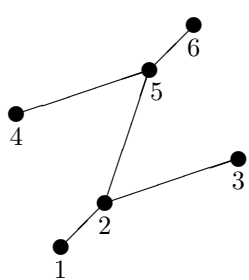
\begin{figure}[ht]
\centering
    \subfigure[permutation $412563$]{
\setlength{\unitlength}{2.8pt}
        \begin{picture}(36,34)(0,0)
        \multiput(0,0)(6,0){7}{\line(0,1){36}}
        \multiput(0,0)(0,6){7}{\line(1,0){36}}
        \put(3, 21){\circle*{2}}
        \put(9,  3){\circle*{2}}
        \put(15, 9){\circle*{2}}
        \put(21,27){\circle*{2}}
        \put(27,33){\circle*{2}}
        \put(33,15){\circle*{2}}
        \end{picture} 
    }
\qquad \qquad \qquad 
    \subfigure[the poset $P_{412563}$]{
\setlength{\unitlength}{2.8pt}
        \begin{picture}(32,34)(3,0)
        \put(9,  3){\line(1,1){6}}
        \put(15, 9){\line(3,1){18}}
        \put(15, 9){\line(1,3){6}}
        \put(3, 21){\line(3,1){18}}
        \put(21,27){\line(1,1){6}}
        \put(9,  3){\circle*{2}}
        \put(15, 9){\circle*{2}}
        \put(3, 21){\circle*{2}}
        \put(21,27){\circle*{2}}
        \put(27,33){\circle*{2}}
        \put(33,15){\circle*{2}}
        \end{picture}
    }
\caption{Permutation posets} 
\label{fig:permposet} 
\end{figure}
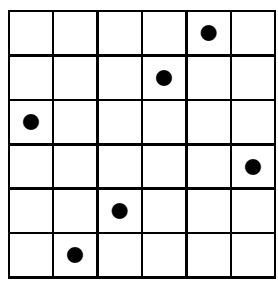
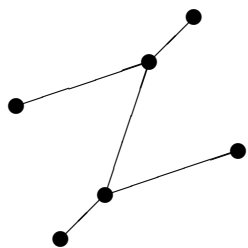


\begin{figure}[ht]
\centering
\setlength{\unitlength}{2.8pt}
\subfigure[$P(\sigma)$]{
        \begin{picture}(36,48)(0,-18)
        \put(9,  3){\line(1,1){6}}
        \put(15, 9){\line(3,1){18}}
        \put(15, 9){\line(1,3){6}}
        \put(3, 21){\line(3,1){18}}
        \put(21,27){\line(1,1){6}}
        \put(3, 21){\circle*{2}}
        \put(9,  3){\circle*{2}}
        \put(15, 9){\circle*{2}}
        \put(21,27){\circle*{2}}
        \put(27,33){\circle*{2}}
        \put(33,15){\circle*{2}}
        \put(9,  0){\makebox(0,0){$1$}} 
        \put(15, 6){\makebox(0,0){$2$}} 
        \put(33,12){\makebox(0,0){$3$}} 
        \put(3, 18){\makebox(0,0){$4$}} 
        \put(22,24){\makebox(0,0){$5$}} 
        \put(27,30){\makebox(0,0){$6$}} 
        \put(6,-18){\line(1,0){12}} 
        \multiput(6,-18)(6,0){3}{\line(0,1){12}} 
        \multiput(6,-12)(0,6){2}{\line(1,0){24}} 
        \multiput(24,-12)(6,0){2}{\line(0,1){6}} 
        \put(9, -9){\makebox(0,0){$1$}} 
        \put(15,-9){\makebox(0,0){$2$}} 
        \put(21,-9){\makebox(0,0){$3$}} 
        \put(9,-15){\makebox(0,0){$4$}} 
        \put(15,-15){\makebox(0,0){$5$}} 
        \put(27,-9){\makebox(0,0){$6$}} 
        \end{picture} 
    }
\qquad \qquad \qquad 
\subfigure[$Q(\sigma)$]{
        \begin{picture}(36,48)(0,-18)
        \put(9,  3){\line(1,1){6}}
        \put(15, 9){\line(3,1){18}}
        \put(15, 9){\line(1,3){6}}
        \put(3, 21){\line(3,1){18}}
        \put(21,27){\line(1,1){6}}
        \put(3, 21){\circle*{2}}
        \put(9,  3){\circle*{2}}
        \put(15, 9){\circle*{2}}
        \put(21,27){\circle*{2}}
        \put(27,33){\circle*{2}}
        \put(33,15){\circle*{2}}
        \put(1, 23){\makebox(0,0){$1$}} 
        \put(7,  5){\makebox(0,0){$2$}} 
        \put(13, 11){\makebox(0,0){$3$}} 
        \put(19,29){\makebox(0,0){$4$}} 
        \put(25,35){\makebox(0,0){$5$}} 
        \put(31,17){\makebox(0,0){$6$}} 
        \put(6,-18){\line(1,0){12}} 
        \multiput(6,-18)(6,0){3}{\line(0,1){12}} 
        \multiput(6,-12)(0,6){2}{\line(1,0){24}} 
        \multiput(24,-12)(6,0){2}{\line(0,1){6}} 
        \put(9, -9){\makebox(0,0){$1$}} 
        \put(9,-15){\makebox(0,0){$2$}} 
        \put(15,-9){\makebox(0,0){$3$}} 
        \put(21,-9){\makebox(0,0){$4$}} 
        \put(27,-9){\makebox(0,0){$5$}} 
        \put(15,-15){\makebox(0,0){$6$}} 
        \end{picture} 
    }
\caption{Tableaux $P(\sigma)$ and $Q(\sigma)$ for $\sigma=412563$} 
\label{fig:PQ} 
\end{figure}

Each permutation poset $P_\sigma$ has two distinguished linear
extensions, obtained by linearly ordering its elements 
$(i,j)=(i,\sigma(i))$
according to the value of the coordinate~$j$ (resp.~$i$), 
as illustrated in Figure~\ref{fig:PQ}. 
The standard tableaux associated with these two linear extensions are
denoted by $P(\sigma)$ and $Q(\sigma)$, respectively. 
The map $\sigma\mapsto (P(\sigma), Q(\sigma))$ is the celebrated 
\emph{Robinson-Schensted correspondence}. 
One of the most striking features of this correspondence is that it is
actually a \emph{bijection} between permutations of~$[n]$ and pairs
$(P,Q)$ of standard Young tableaux that have the same shape consisting
of $n$ boxes. 

We will next explain why the conventional description of this
correspondence, due to C.~Schensted~\cite{schensted}
(cf.~\cite[Section~7.11]{stanley2}), is equivalent 
to the one we just gave,
as first observed and proved by 
C.~Greene (see Theorem~\ref{th:greene74} below). 

Fix a permutation poset $P_\sigma\,$, and consider its order
ideals $P_\sigma(i,j)$ defined by
\begin{equation}
\label{eq:P(i,j)}
P_\sigma(i,j) = ([i]\times[j])\cap P_\sigma \,, 
\end{equation}
for $i,j\in \{0,1,\dots,n\}$.
Thus $P_\sigma(i,j)$ consists of the points $(k,\sigma(k))$
located (weakly) southwest of~$(i,j)$. 
The shapes of these order ideals are denoted by
\[
\dc_{ij}=\dc(P_\sigma(i,j)) \,.
\]
The two-dimensional array $(\dc_{ij})$ is called the \emph{growth
  diagram} for the permutation~$\sigma$. 
An example of a growth diagram is given in Figure~\ref{fig:growP}.

\begin{figure}[ht]
\centering 
\setlength{\unitlength}{2.7pt}
\begin{picture}(83,83)(-11,-11)
\multiput(-3,-4)(0,12){7}{\line(1,0){72}}
\multiput(-3,-4)(12,0){7}{\line(0,1){72}}
\put(-10,-4){\makebox(0,0){$0$}}
\put(-10, 8){\makebox(0,0){$1$}}
\put(-10,20){\makebox(0,0){$2$}}
\put(-10,32){\makebox(0,0){$3$}}
\put(-10,44){\makebox(0,0){$4$}}
\put(-10,56){\makebox(0,0){$5$}}
\put(-10,68){\makebox(0,0){$6$}}
\put(-3,-12){\makebox(0,0){$0$}}
\put( 9,-12){\makebox(0,0){$1$}}
\put(21,-12){\makebox(0,0){$2$}}
\put(33,-12){\makebox(0,0){$3$}}
\put(45,-12){\makebox(0,0){$4$}}
\put(57,-12){\makebox(0,0){$5$}}
\put(69,-12){\makebox(0,0){$6$}}
\multiput(16.5, 2)(12, 12){2}{\makebox(0,0){\circle*{2.5}}}
\multiput(40.5,50)(12, 12){2}{\makebox(0,0){\circle*{2.5}}}
\multiput( 4.5,38)(60,-12){2}{\makebox(0,0){\circle*{2.5}}}
\multiput(-4.5,-6)(0,12){7}{\makebox(0,0){$\phi$}} 
\multiput( 7.5,-6)(12,0){6}{\makebox(0,0){$\phi$}} 
\multiput( 7.5, 6)(0,12){3}{\makebox(0,0){$\phi$}}
\multiput( 7.5,42.5)(0,12){3}{\makebox(0,0){\begin{picture}(1,1)(0,0)
                                       \multiput(0,0)(0,1){2}{\line(1,0){1}}
                                       \multiput(0,0)(1,0){2}{\line(0,1){1}}
                                       \end{picture}}}
\multiput(19.5,6.5)(0,12){3}{\makebox(0,0){\begin{picture}(1,1)(0,0)
                                       \multiput(0,0)(0,1){2}{\line(1,0){1}}
                                       \multiput(0,0)(1,0){2}{\line(0,1){1}}
                                       \end{picture}}}
\multiput(31.5,6.5)(12,0){4}{\makebox(0,0){\begin{picture}(1,1)(0,0)
                                       \multiput(0,0)(0,1){2}{\line(1,0){1}}
                                       \multiput(0,0)(1,0){2}{\line(0,1){1}}
                                       \end{picture}}}
\multiput(19.5,42)(0,12){3}{\makebox(0,0){\begin{picture}(1,2)(0,0)
                                       \multiput(0,0)(0,1){3}{\line(1,0){1}}
                                       \multiput(0,0)(1,0){2}{\line(0,1){2}}
                                       \end{picture}}}
\multiput(31,30.5)(12,0){3}{\makebox(0,0){\begin{picture}(2,1)(0,0)
                                       \multiput(0,0)(0,1){2}{\line(1,0){2}}
                                       \multiput(0,0)(1,0){3}{\line(0,1){1}}
                                       \end{picture}}}
\multiput(31,18.5)(12,0){4}{\makebox(0,0){\begin{picture}(2,1)(0,0)
                                       \multiput(0,0)(0,1){2}{\line(1,0){2}}
                                       \multiput(0,0)(1,0){3}{\line(0,1){1}}
                                       \end{picture}}}
\put(66.5,30.5){\makebox(0,0){\begin{picture}(3,1)(0,0)
                          \multiput(0,0)(0,1){2}{\line(1,0){3}}
                          \multiput(0,0)(1,0){4}{\line(0,1){1}}
                          \end{picture}}}
\multiput(31,42)(0,12){3}{\makebox(0,0){\begin{picture}(2,2)(0,0)
                                       \put(0,0){\line(1,0){1}}
                                       \put(2,1){\line(0,1){1}}
                                       \multiput(0,1)(0,1){2}{\line(1,0){2}}
                                       \multiput(0,0)(1,0){2}{\line(0,1){2}}
                                       \end{picture}}}
\multiput(43,42)(12,0){2}{\makebox(0,0){\begin{picture}(2,2)(0,0)
                                       \put(0,0){\line(1,0){1}}
                                       \put(2,1){\line(0,1){1}}
                                       \multiput(0,1)(0,1){2}{\line(1,0){2}}
                                       \multiput(0,0)(1,0){2}{\line(0,1){2}}
                                       \end{picture}}}
\multiput(54.5,54)(12,-12){2}{\makebox(0,0){\begin{picture}(3,2)(0,0)
                                       \put(0,0){\line(1,0){1}}
                                       \multiput(2,1)(1,0){2}{\line(0,1){1}}
                                       \multiput(0,1)(0,1){2}{\line(1,0){3}}
                                       \multiput(0,0)(1,0){2}{\line(0,1){2}}
                                       \end{picture}}}
\multiput(42.5,54)(0,12){2}{\makebox(0,0){\begin{picture}(3,2)(0,0)
                                       \put(0,0){\line(1,0){1}}
                                       \multiput(2,1)(1,0){2}{\line(0,1){1}}
                                       \multiput(0,1)(0,1){2}{\line(1,0){3}}
                                       \multiput(0,0)(1,0){2}{\line(0,1){2}}
                                       \end{picture}}}
\put(54,66){\makebox(0,0){\begin{picture}(4,2)(0,0)
                          \put(0,0){\line(1,0){1}}
                          \multiput(2,1)(1,0){3}{\line(0,1){1}}
                          \multiput(0,1)(0,1){2}{\line(1,0){4}}
                          \multiput(0,0)(1,0){2}{\line(0,1){2}}
                          \end{picture}}}
\put(66.5,54){\makebox(0,0){\begin{picture}(3,2)(0,0)
                          \put(0,0){\line(1,0){2}}
                          \put(3,1){\line(0,1){1}}
                          \multiput(0,1)(0,1){2}{\line(1,0){3}}
                          \multiput(0,0)(1,0){3}{\line(0,1){2}}
                          \end{picture}}}
\put(66,66){\makebox(0,0){\begin{picture}(4,2)(0,0)
                          \put(0,0){\line(1,0){2}}
                          \multiput(3,1)(1,0){2}{\line(0,1){1}}
                          \multiput(0,1)(0,1){2}{\line(1,0){4}}
                          \multiput(0,0)(1,0){3}{\line(0,1){2}}
                          \end{picture}}}
\end{picture}
\caption{The growth diagram for $\sigma=412563$}
\label{fig:growP}
\end{figure}
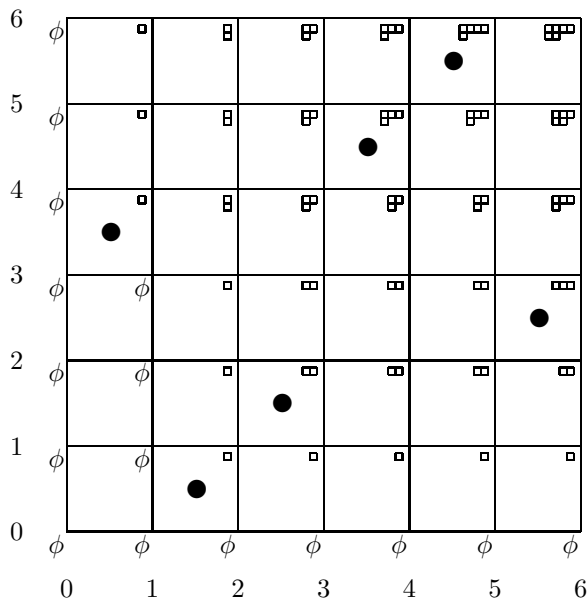

Remarkably, the shapes that make up any growth diagram satisfy a very
simple local rule that provides a recursive algorithm for
computing~$\dc(P_\sigma)$. 
Specifically, let us consider an arbitrary $2\times 2$ submatrix
\begin{equation}
\label{eq:2by2submatrix}
\begin{matrix} 
 \dc_{i-1, j} & \dc_{ij} \\[.2in]
 \dc_{i-1, j-1}   & \dc_{i, j-1} 
\end{matrix}
\end{equation}
of the growth diagram for a permutation~$\sigma$. 
The following theorem shows that the shape $\dc_{ij}$ is uniquely
determined by the shapes $\dc_{i-1, j-1}$, $\dc_{i, j-1}$, 
and $\dc_{i-1, j}\,$, together with knowing whether $\sigma(i)=j$ or
not (i.e., whether $(i,j)\in P_\sigma$ or not).

\begin{theorem} \label{thm:growP} \cite{fomin3} 
\begin{enumerate}
\item 
If $\dc_{i,j-1} \neq \dc_{i-1,j}$,  
then $\dc_{ij}=\dc_{i,j-1} \cup \dc_{i-1,j}\,$. 
\item
If $\dc_{i,j-1}=\dc_{i-1,j}=\dc_{i-1,j-1}$  
and $\sigma(i)\neq j$, then $\dc_{ij}=\dc_{i-1,j-1}\,$. 
\item 
If $\dc_{i,j-1}=\dc_{i-1,j}=\dc_{i-1,j-1}$  
and $\sigma(i)=j$, then $\dc_{ij}$ is obtained 
by adding a box to the first row of $\dc_{i-1,j-1}\,$. 
\item 
If $\dc_{i,j-1}=\dc_{i-1,j}\neq \dc_{i-1,j-1}\,$,
then $\dc_{ij}$ is obtained by adding a box to the row 
immediately below the box $\dc_{i-1,j}-\dc_{i-1,j-1}\,$. 
\end{enumerate}
\end{theorem} 

\begin{proof}
Part~1 follows from the Monotonicity Theorem.

Assume $\dc_{i,j-1}=\dc_{i-1,j}=\dc_{i-1,j-1}\,$. 
Then the $i$th column (resp.\ $j$th row) does not contain elements of
$P_\sigma$ below (resp.\ to the left) of $(i, j)$. 
If, in addition, $(i, j)\notin P_\sigma\,$, 
then $P_\sigma(i-1,j-1)=P_\sigma(i,j)$ and therefore 
$\dc_{ij}=\dc_{i-1, j-1}\,$, proving Part~2. 
If, on the other hand, $(i, j)\in P_\sigma\,$ 
(i.e., $\sigma(i)=j$), then $(i, j)$
is greater than all elements strictly below and to the left of it, 
so any chain in $P_\sigma(i-1,j-1)$ is extended by~$p$. 
Hence the maximal length, $c=c_1$, of a chain is increased by~1. 
The first row of $\dc_{ij}$ then contains one more box 
than the first row of $\dc_{i-1, j-1}$, proving Part~3. 

To prove Part~4, assume $\dc_{i,j-1}=\dc_{i-1,j}\neq \dc_{i-1,j-1}\,$. 
Then $P_\sigma$ contains an element $p_1$ 
strictly to the left of $(i, j)$, 
as well as an element $p_2$ strictly below $(i, j)$. 
Both $p_1$ and $p_2$ are maximal elements of $P_\sigma(i,j)\,$, 
and $\dc (P_\sigma(i, j)-\{p_1\})=\dc_{i,j-1}
=\dc_{i-1,j}=\dc (P_\sigma(i, j)-\{p_2\})$. 
Let the boxes $A$ and $B$ be defined by 
$\{B\}=\dc_{ij}-\dc_{i-1,j}$ and 
$\{A\}=\dc_{i-1,j}-\dc_{i-1, j-1}\,$. 
Theorem~\ref{thm:fominG} implies that $B$ lies weakly to 
the left of $A$ (see Figure~\ref{fig:proof1}a). 
Now consider the poset $P'_\sigma(i,j)$ on the same ground set as 
$P_\sigma(i,j)$, the difference being that 
$(k,\sigma(k))\leq (l,\sigma(l))$ in $P'_\sigma(i,j)$ if and only if 
$k\geq l$ and $\sigma(k)\leq \sigma(l)$.  
The chains of $P'_\sigma(i,j)$ are the antichains of~$P_\sigma(i,j)$. 
Hence the shape 
$\dc_{ij} '=\dc(P'_{\sigma}(i,j))$ 
is the transpose of the shape~$\dc_{ij}\,$. 
Notice that in $P'_\sigma(i,j)$, 
$p_1$ is maximal while $p_2$ is minimal. 
Theorem~\ref{thm:fominG} then implies that $B$ lies in either 
the same \emph{row} as $A$, 
or in the row immediately below it (see Figure~\ref{fig:proof1}b).  
We conclude that $B$ must lie one row below~$A$ 
(see Figure~\ref{fig:proof1}c), as desired.  
\end{proof}

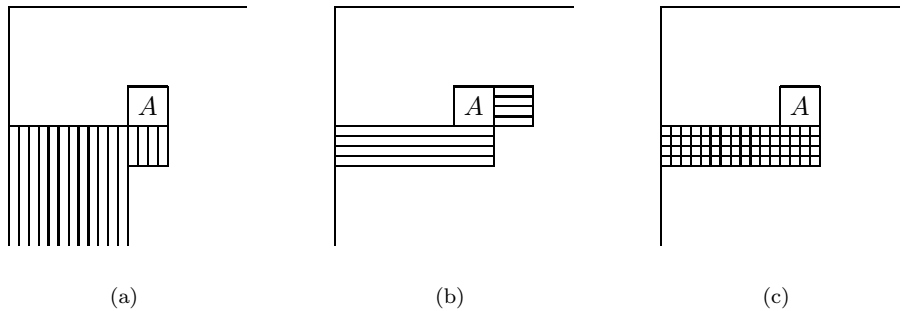
\begin{figure}[ht]
\centering
\setlength{\unitlength}{2.5pt}
\mbox{
   \subfigure[]{ 
        \begin{picture}(36,36)(0,0)
        \put(0, 0){\line(0,1){36}}
        \put(0,36){\line(1,0){36}}
        \put(0,18){\line(1,0){24}}
        \put(18,0){\line(0,1){12}}
        \put(18,18){\line(0,1){6}}
        \put(24,12){\line(0,1){12}}
        \multiput(18,12)(0,12){2}{\line(1,0){6}}
        \put(21,21){\makebox(0,0){$A$}} 
        \multiput(0, 0)(1.5,0){13}{\line(0,1){18}} 
        \multiput(18,12)(1.5,0){5}{\line(0,1){6}} 
        \end{picture} 
    }
\qquad \
    \subfigure[]{
        \begin{picture}(36,36)(0,0)
        \put(0, 0){\line(0,1){36}}
        \put(0,36){\line(1,0){36}}
        \put(24,12){\line(0,1){12}} 
        \put(30,18){\line(0,1){6}} 
        \multiput(0 ,12)(0,6){2}{\line(1,0){24}}
        \multiput(18,18)(0,6){2}{\line(1,0){12}}
        \put(18,18){\line(0,1){6}}
        \put(21,21){\makebox(0,0){$A$}} 
        \multiput(0 ,12)(0,1.5){5}{\line(1,0){24}} 
        \multiput(24,18)(0,1.5){5}{\line(1,0){6}} 
        \end{picture}
    }
\qquad \
    \subfigure[]{
        \begin{picture}(36,36)(0,0)
        \put(0, 0){\line(0,1){36}}
        \put(0,36){\line(1,0){36}}
        \multiput(18,18)(6,0){2}{\line(0,1){6}}
        \multiput(18,18)(0,6){2}{\line(1,0){6}} 
        \put(21,21){\makebox(0,0){$A$}} 
        \multiput(0,12)(0,1.5){5}{\line(1,0){24}} 
        \multiput(0,12)(1.5,0){17}{\line(0,1){6}} 
        \end{picture}
   }
} 
\caption{Allowable locations of~$B$} 
\label{fig:proof1} 
\end{figure}


The growth rules described in Theorem~\ref{thm:growP} 
can be used recursively to compute the shape 
$\dc(P_\sigma)=|P(\sigma)|=|Q(\sigma)|$, beginning by putting the
empty shapes at the southwest border of the growth diagram, and
expanding northeast with the help of recursion. 
This parallel algorithm has various sequential versions;
let us choose the one where the shapes $\dc_{ij}$ are computed column
by column (left to right; and bottom-up within each column).
Let $P_i$ denote the (non-standard) tableau that encodes the $i$th
column of the growth diagram, for $i=0,1,\dots,n$;
more precisely, this tableau has entry $j$ in a box $B$ provided 
$\{B\}=\dc_{ij}-\dc_{i,j-1}\,$. 
Figure~\ref{fig:RS} shows the tableaux $P_i$ for our
running example $\sigma=(412563)$. 
(These tableaux are obtained by encoding the columns of the growth
diagram in Figure~\ref{fig:growP}.) 

By the nature of the recursion process, the tableau $P_i$ is
completely determined by the previous tableau $P_{i-1}$ together with
the entry $\sigma(i)$ of the permutation~$\sigma$. 
The rule for computing $P_i$ from $P_{i-1}$ and $\sigma(i)$
can be reformulated entirely in the language of tableaux;
one then arrives at the familiar ``insertion'' step of Schensted's
algorithm~\cite{schensted}. (We leave this verification to the reader.)
Thus the tableau $P_n=P(\sigma)$ is indeed Schensted's $P$-tableau
(sometimes called the ``insertion tableau'') for~$\sigma$. 
The growth of the shapes $\dc_{in}$ of the tableaux $P_i$ is recorded
by $Q(\sigma)$, which is therefore the $Q$-tableau (or the
``recording tableau'') of Schensted's original construction. 

\begin{figure}[ht]
\centering 
\setlength{\unitlength}{2pt}
\begin{picture}(160,50)(0,0)
\put( 0,46){\makebox(0,0){$P_0$}} 
\put(16,46){\makebox(0,0){$P_1$}} 
\put(33,46){\makebox(0,0){$P_2$}} 
\put(53,46){\makebox(0,0){$P_3$}} 
\put(80,46){\makebox(0,0){$P_4$}} 
\put(112,46){\makebox(0,0){$P_5$}} 
\put(147,46){\makebox(0,0){$P_6=P(\sigma)$}} 
\put(147,0){\makebox(0,0){$Q(\sigma)$}} 

\put(0,35){\makebox(0,0){$\phi$}}
\put(16,35){\makebox(0,0){$4$}}
\multiput(33,29)(17,0){2}{\makebox(0,0){$4$}} 
\multiput(33,35)(17,0){2}{\makebox(0,0){$1$}} 
\put(56,35){\makebox(0,0){$2$}} 
\multiput(74,29)(29,0){2}{\makebox(0,0){$4$}} 
\multiput(74,35)(29,0){2}{\makebox(0,0){$1$}} 
\multiput(80,35)(29,0){2}{\makebox(0,0){$2$}} 
\multiput(86,35)(29,0){2}{\makebox(0,0){$5$}} 
\multiput(156,17)(35,0){1}{\makebox(0,0){$5$}} 
\multiput(121,35)(35,0){2}{\makebox(0,0){$6$}} 
\multiput(138,35)(0,-18){2}{\makebox(0,0){$1$}} 
\multiput(144,35)(-6,-24){2}{\makebox(0,0){$2$}} 
\multiput(150,35)(-6,-18){2}{\makebox(0,0){$3$}} 
\multiput(138,29)(12,-12){2}{\makebox(0,0){$4$}} 
\put(144,11){\makebox(0,0){$6$}} 
\put(144,29){\makebox(0,0){$5$}}  
\multiput(16,32)(0,18){1}{\makebox(0,0){\begin{picture}(6,12)(0,0) 
                        \multiput(0,6)(0,6){2}{\line(1,0){6}} 
                        \multiput(0,6)(6,0){2}{\line(0,1){6}} 
                        \end{picture}}} 
\multiput(33,32)(0,18){1}{\makebox(0,0){\begin{picture}(6,12)(0,0) 
                        \multiput(0,0)(0,6){3}{\line(1,0){6}} 
                        \multiput(0,0)(6,0){2}{\line(0,1){12}} 
                        \end{picture}}} 
\multiput(53,32)(0,18){1}{\makebox(0,0){\begin{picture}(12,12)(0,0) 
                        \put(0,0){\line(1,0){6}} 
                        \put(12,6){\line(0,1){6}} 
                        \multiput(0,0)(6,0){2}{\line(0,1){12}} 
                        \multiput(0,6)(0,6){2}{\line(1,0){12}} 
                        \end{picture}}} 
\multiput(80,32)(0,18){1}{\makebox(0,0){\begin{picture}(18,12)(0,0) 
                        \put(0,0){\line(1,0){6}} 
                        \multiput(12,6)(6,0){2}{\line(0,1){6}} 
                        \multiput(0,0)(6,0){2}{\line(0,1){12}} 
                        \multiput(0,6)(0,6){2}{\line(1,0){18}} 
                        \end{picture}}} 
\multiput(112,32)(0,18){1}{\makebox(0,0){\begin{picture}(24,12)(0,0) 
                        \put(0,0){\line(1,0){6}} 
                        \multiput(12,6)(6,0){3}{\line(0,1){6}} 
                        \multiput(0,0)(6,0){2}{\line(0,1){12}} 
                        \multiput(0,6)(0,6){2}{\line(1,0){24}} 
                        \end{picture}}} 
\multiput(147,14)(0,18){2}{\makebox(0,0){\begin{picture}(24,12)(0,0) 
                        \put(0,0){\line(1,0){12}} 
                        \multiput(18,6)(6,0){2}{\line(0,1){6}} 
                        \multiput(0,0)(6,0){3}{\line(0,1){12}} 
                        \multiput(0,6)(0,6){2}{\line(1,0){24}} 
                        \end{picture}}} 
\end{picture} 
\caption{The tableaux $P_i$ for $\sigma=412563$}
\label{fig:RS}
\end{figure}
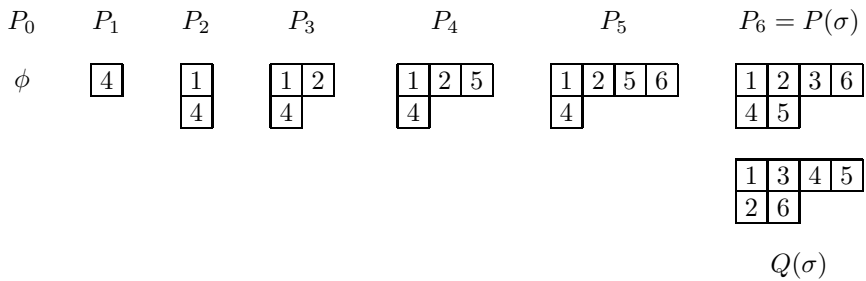 

We thus obtained the following fundamental result of 
C.~Greene. 

\begin{theorem}
\label{th:greene74}
\cite{greene2}
The common shape of the two tableaux 
associated to a given permutation $\sigma$ by the Robinson-Schensted
correspondence (as described by Schensted~\cite{schensted}) 
is exactly the shape $\dc(P_\sigma)$ for the permutation
poset~$P_\sigma\,$. 
\end{theorem}

Several properties of the Robinson-Schensted algorithm, 
which are hard to obtain from the usual ``bumping'' description, 
are easily deduced from the growth diagram approach. 
Here are two examples. 

\begin{corollary} 
\cite{schensted, schutzenberger}
\label{cor:inv}
Inverting a permutation interchanges the two tableaux associated to
it: $P(\sigma^{-1})=Q(\sigma )$, $Q(\sigma^{-1})=P(\sigma )$. 
\end{corollary}

\proof 
Transpose the growth diagram of $\sigma$ in the diagonal that connects
its southwest and northeast corners. 
Then $\sigma$ becomes $\sigma^{-1}$, 
while $P(\sigma)$ and $Q(\sigma)$ are interchanged. 
\endproof

\begin{corollary} \label{cor:bij}
The map $\sigma\mapsto (P(\sigma ),Q(\sigma ))$ 
is a bijection between permutations of $[n]$, on one hand, 
and pairs of standard Young tableaux of the same shape consisting of $n$ 
boxes, on the other. 
\end{corollary}

\proof 
(Sketch)
It is straightforward to verify, using Theorem~\ref{thm:growP}, that
for any $2\times 2$ submatrix (\ref{eq:2by2submatrix}) in the growth
diagram, the shape $\dc_{i-1,j-1}$ is uniquely determined by the three
shapes $\dc_{i-1,j}$, $\dc_{i,j-1}$, and $\dc_{ij}$, and furthermore
these three shapes determine whether $\sigma(i)=j$ or not. 
Thus the whole growth diagram can be reconstructed recursively, 
beginning at the upper-right boundary
(i.e., using $P(\sigma)$ and $Q(\sigma)$ as inputs);
along the way,
we will recover~$\sigma$,
as desired. 
\endproof


\end{section}

\begin{section}{The Sch\"utzenberger involution}  
\label{sec:schutzenberger} 

Recall that the $Q$-tableaux $Q(\sigma )$ is obtained by ``growing'' 
the permutation poset~$P_\sigma$ (thus the corresponding shape) along
the ``left-to-right'' linear extension, as shown in
Figure~\ref{fig:PQ}b (or in Figure~\ref{fig:Q'}a below). 
Alternatively, we could have grown the poset from right to left,
by consecutively adding the elements labelled $n,n-1,n-2,\dots$, 
in this order. 
The corresponding standard Young tableaux $Q'(\sigma)$ has of course
the same shape as $Q(\sigma)$; see Figure~\ref{fig:Q'}b. 
Notice that $Q'(\sigma)=Q(\sigma')$,
where $\sigma'$ is the permutation defined by 
$\sigma '(i)=n\!+\!1\!-\!\sigma(n\!+\!1\!-i)$. 
(In other words, $\sigma'$ is obtained from $\sigma$ by 180
degrees rotation.) See Figure~\ref{fig:Q'}c. 

Remarkably, the tableaux $Q'(\sigma )$ can be computed from
$Q(\sigma)$ alone, without knowing $\sigma$ itself. 
The corresponding construction is the famous 
\emph{Sch\"utzenberger involution}, as we explain below. 
(This presentation follows \cite{fomin3}; cf.\ also \cite{fomin6,
  leeuwen1}.)

\begin{figure}[ht]
\centering
\setlength{\unitlength}{2.8pt} 
\subfigure[$Q(126453)$]{
        \begin{picture}(36,58)(0,-24)
        \put(3,3){\line(1,1){24}}
        \put(9,9){\line(4,1){24}}
        \put(9,9){\line(1,4){6}}
        \multiput(3,3)(6,6){2}{\circle*{2}} 
        \multiput(21,21)(6,6){2}{\circle*{2}} 
        \put(15,33){\circle*{2}} 
        \put(33,15){\circle*{2}} 
        \put(1,  5){\makebox(0,0){$1$}} 
        \put(7, 11){\makebox(0,0){$2$}} 
        \put(13,35){\makebox(0,0){$3$}} 
        \put(19,23){\makebox(0,0){$4$}} 
        \put(25,29){\makebox(0,0){$5$}} 
        \put(31,17){\makebox(0,0){$6$}} 
        \multiput(6,-24)(0,6){2}{\line(1,0){6}} 
        \multiput(6,-24)(6,0){2}{\line(0,1){18}} 
        \multiput(6,-12)(0,6){2}{\line(1,0){24}} 
        \multiput(18,-12)(6,0){3}{\line(0,1){6}} 
        \put(9, -9){\makebox(0,0){$1$}} 
        \put(9,-15){\makebox(0,0){$4$}} 
        \put(15,-9){\makebox(0,0){$2$}} 
        \put(9,-21){\makebox(0,0){$6$}} 
        \put(21,-9){\makebox(0,0){$3$}} 
        \put(27,-9){\makebox(0,0){$5$}} 
        \end{picture} 
    }
\quad 
\subfigure[$Q'(126453)$]{
        \begin{picture}(36,58)(0,-24)
        \put(3,3){\line(1,1){24}}
        \put(9,9){\line(4,1){24}}
        \put(9,9){\line(1,4){6}}
        \multiput(3,3)(6,6){2}{\circle*{2}} 
        \multiput(21,21)(6,6){2}{\circle*{2}} 
        \put(15,33){\circle*{2}} 
        \put(33,15){\circle*{2}} 
        \put(1,  5){\makebox(0,0){$6$}} 
        \put(7, 11){\makebox(0,0){$5$}} 
        \put(13,35){\makebox(0,0){$4$}} 
        \put(19,23){\makebox(0,0){$3$}} 
        \put(25,29){\makebox(0,0){$2$}} 
        \put(31,17){\makebox(0,0){$1$}} 
        \multiput(6,-24)(0,6){2}{\line(1,0){6}} 
        \multiput(6,-24)(6,0){2}{\line(0,1){18}} 
        \multiput(6,-12)(0,6){2}{\line(1,0){24}} 
        \multiput(18,-12)(6,0){3}{\line(0,1){6}} 
        \put(9, -9){\makebox(0,0){$1$}} 
        \put(9,-15){\makebox(0,0){$2$}} 
        \put(15,-9){\makebox(0,0){$3$}} 
        \put(9,-21){\makebox(0,0){$4$}} 
        \put(21,-9){\makebox(0,0){$5$}} 
        \put(27,-9){\makebox(0,0){$6$}} 
        \end{picture} 
    }
\quad 
\subfigure[$Q(423156)$]{
        \begin{picture}(36,58)(0,-24)
        \put(9, 9){\line(1,1){24}}
        \put(3,21){\line(4,1){24}}
        \put(21,3){\line(1,4){6}} 
        \multiput(3,21)(12,-6){2}{\circle*{2}} 
        \multiput(9, 9)(12,-6){2}{\circle*{2}} 
        \multiput(27,27)(6,6){2}{\circle*{2}} 
        \put(1, 23){\makebox(0,0){$1$}} 
        \put(7, 11){\makebox(0,0){$2$}} 
        \put(13,17){\makebox(0,0){$3$}} 
        \put(19, 5){\makebox(0,0){$4$}} 
        \put(25,29){\makebox(0,0){$5$}} 
        \put(31,35){\makebox(0,0){$6$}} 
        \multiput(6,-24)(0,6){2}{\line(1,0){6}} 
        \multiput(6,-24)(6,0){2}{\line(0,1){18}} 
        \multiput(6,-12)(0,6){2}{\line(1,0){24}} 
        \multiput(18,-12)(6,0){3}{\line(0,1){6}} 
        \put(9, -9){\makebox(0,0){$1$}} 
        \put(9,-15){\makebox(0,0){$2$}} 
        \put(15,-9){\makebox(0,0){$3$}} 
        \put(9,-21){\makebox(0,0){$4$}} 
        \put(21,-9){\makebox(0,0){$5$}} 
        \put(27,-9){\makebox(0,0){$6$}} 
        \end{picture} 
    }
\caption{Tableaux $Q(\sigma)$, $Q'(\sigma)$, and $Q(\sigma ')$} 
\label{fig:Q'} 
\end{figure}
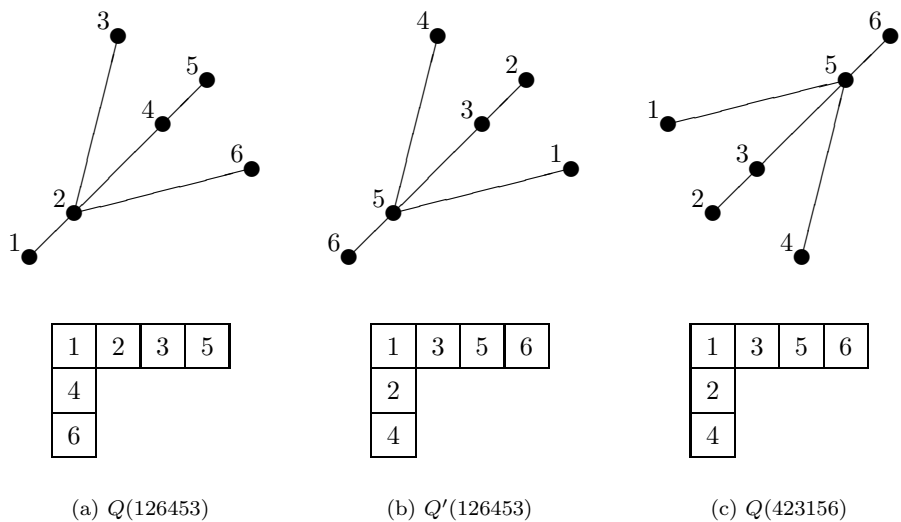

For $1\leq i\leq j\leq n$, let us define the subposet $P_\sigma[i;j]$
of $P_\sigma$ by 
\[
P_\sigma[i;j]= 
\{(k,\sigma (k))\,:\, i\leq k\leq j\} \,. 
\]
The corresponding shapes will be denoted by 
$\dc_{[i;j]}=\dc(P_\sigma[i;j])$ 
and placed in a triangular array, as shown in
Figure~\ref{fig:collat}. 
Note that the upper-left side of this array is encoded by the tableau
$Q(\sigma)$, while the upper-right side corresponds to~$Q'(\sigma)$. 

As in Section~\ref{sec:schensted}, 
Theorem~\ref{thm:fominG} can be used to 
obtain local rules of growth in this array. 
Let us fix $1< i\leq j<n$, and consider the following four shapes: 
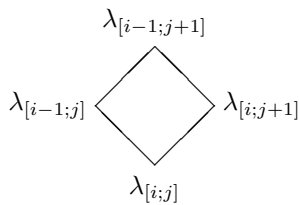
\begin{align*} 
\setlength{\unitlength}{4.5pt} 
\begin{picture}(18,17)(6,2) 
\centering 
\multiput(10,10)(5,-5){2}{\line(1,1){5}}
\multiput(15,5)(5,5){2}{\line(-1,1){5}}
\put(15,3){\makebox(0,0){$\dc_{[i;j]}$}}
\put(6,10){\makebox(0,0){$\dc_{[i-1;j]}$}}
\put(24,10){\makebox(0,0){$\dc_{[i;j+1]}$}}
\put(15,17){\makebox(0,0){$\dc_{[i-1;j+1]}$}}
\end{picture} 
\end{align*}

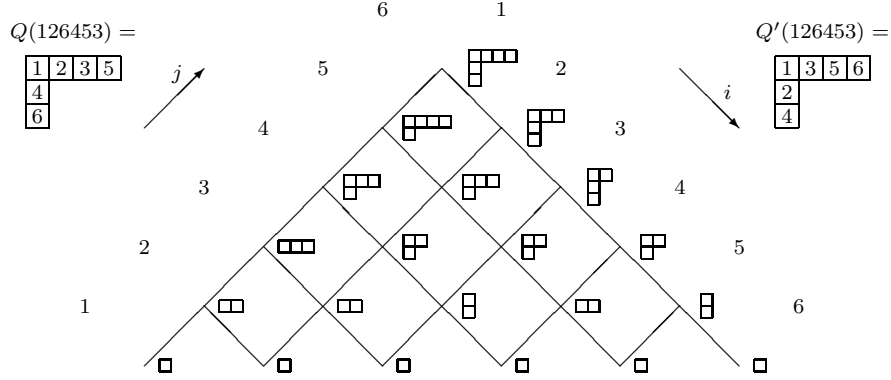
\begin{figure}[ht] 
\centering 
\setlength{\unitlength}{4.5pt} 
\footnotesize 
\begin{picture}(52,43)(-1,5)
\put(0,10){\line(1,1){25}} 
\put(10,10){\line(1,1){20}} 
\put(20,10){\line(1,1){15}} 
\put(30,10){\line(1,1){10}} 
\put(40,10){\line(1,1){5}} 
\put(10,10){\line(-1,1){5}} 
\put(20,10){\line(-1,1){10}} 
\put(30,10){\line(-1,1){15}} 
\put(40,10){\line(-1,1){20}} 
\put(50,10){\line(-1,1){25}} 
\put(45,35){\vector(1,-1){5}} 
\put(0,30){\vector(1,1){5}} 
\put(2.75,34.5){\makebox(0,0){$j$}}
\put(49,33){\makebox(0,0){$i$}}
\put(-5,15){\makebox(0,0){$1$}}
\put(0 ,20){\makebox(0,0){$2$}}
\put( 5,25){\makebox(0,0){$3$}}
\put(10,30){\makebox(0,0){$4$}}
\put(15,35){\makebox(0,0){$5$}}
\put(20,40){\makebox(0,0){$6$}}
\put(30,40){\makebox(0,0){$1$}}
\put(35,35){\makebox(0,0){$2$}}
\put(40,30){\makebox(0,0){$3$}}
\put(45,25){\makebox(0,0){$4$}}
\put(50,20){\makebox(0,0){$5$}}
\put(55,15){\makebox(0,0){$6$}}
\multiput(1.75 ,10)(10,0){6}{\makebox(0,0){\begin{picture}(1,1)(0,0)
                                \multiput(0,0)(0,1){2}{\line(1,0){1}}
                                \multiput(0,0)(1,0){2}{\line(0,1){1}}
                                \end{picture}}}
\multiput(7.25 ,15)(10,0){2}{\makebox(0,0){\begin{picture}(2,1)(0,0) 
                                \multiput(0,0)(0,1){2}{\line(1,0){2}} 
                                \multiput(0,0)(1,0){3}{\line(0,1){1}} 
                                \end{picture}}} 
\multiput(27.25,15)(20,0){2}{\makebox(0,0){\begin{picture}(1,2)(0,0) 
                                \multiput(0,0)(0,1){3}{\line(1,0){1}} 
                                \multiput(0,0)(1,0){2}{\line(0,1){2}} 
                                \end{picture}}} 
\put(37.25,15){\makebox(0,0){\begin{picture}(2,1)(0,0) 
                                \multiput(0,0)(0,1){2}{\line(1,0){2}} 
                                \multiput(0,0)(1,0){3}{\line(0,1){1}} 
                                \end{picture}}} 
\put(12.75,20){\makebox(0,0){\begin{picture}(3,1)(0,0) 
                                \multiput(0,0)(0,1){2}{\line(1,0){3}} 
                                \multiput(0,0)(1,0){4}{\line(0,1){1}} 
                                \end{picture}}} 
\multiput(22.75,20)(10,0){3}{\makebox(0,0){\begin{picture}(2,2)(0,0) 
                                \put(0,0){\line(1,0){1}} 
                                \put(2,1){\line(0,1){1}} 
                                \multiput(0,1)(0,1){2}{\line(1,0){2}} 
                                \multiput(0,0)(1,0){2}{\line(0,1){2}} 
                                \end{picture}}} 
\multiput(18.25,25)(10,0){2}{\makebox(0,0){\begin{picture}(3,2)(0,0) 
                                \put(0,0){\line(1,0){1}} 
                                \multiput(2,1)(1,0){2}{\line(0,1){1}} 
                                \multiput(0,1)(0,1){2}{\line(1,0){3}} 
                                \multiput(0,0)(1,0){2}{\line(0,1){2}} 
                                \end{picture}}} 
\put(38.25,25){\makebox(0,0){\begin{picture}(2,3)(0,0) 
                                \put(2,2){\line(0,1){1}} 
                                \multiput(0,2)(0,1){2}{\line(1,0){2}} 
                                \multiput(0,0)(0,1){2}{\line(1,0){1}} 
                                \multiput(0,0)(1,0){2}{\line(0,1){3}} 
                                \end{picture}}} 
\put(23.75,30){\makebox(0,0){\begin{picture}(4,2)(0,0) 
                                \put(0,0){\line(1,0){1}} 
                                \multiput(2,1)(1,0){3}{\line(0,1){1}} 
                                \multiput(0,1)(0,1){2}{\line(1,0){4}} 
                                \multiput(0,0)(1,0){2}{\line(0,1){2}} 
                                \end{picture}}} 
\put(33.75,30){\makebox(0,0){\begin{picture}(3,3)(0,0) 
                                \multiput(0,0)(0,1){2}{\line(1,0){1}} 
                                \multiput(2,2)(1,0){2}{\line(0,1){1}} 
                                \multiput(0,2)(0,1){2}{\line(1,0){3}} 
                                \multiput(0,0)(1,0){2}{\line(0,1){3}} 
                                \end{picture}}} 
\put(29.25,35){\makebox(0,0){\begin{picture}(4,3)(0,0) 
                                \multiput(0,0)(0,1){2}{\line(1,0){1}} 
                                \multiput(2,2)(1,0){3}{\line(0,1){1}} 
                                \multiput(0,2)(0,1){2}{\line(1,0){4}} 
                                \multiput(0,0)(1,0){2}{\line(0,1){3}} 
                                \end{picture}}} 
\put(-6,38){\makebox(0,0){$Q(126453)=$}} 
\put(-6,33){\makebox(0,0){\begin{picture}(8,6)(0,0) 
                                \multiput(0,0)(0,2){2}{\line(1,0){2}} 
                                \multiput(4,4)(2,0){3}{\line(0,1){2}} 
                                \multiput(0,4)(0,2){2}{\line(1,0){8}} 
                                \multiput(0,0)(2,0){2}{\line(0,1){6}} 
                                \put(1,5){\makebox(0,0){$1$}} 
                                \put(3,5){\makebox(0,0){$2$}} 
                                \put(5,5){\makebox(0,0){$3$}} 
                                \put(1,3){\makebox(0,0){$4$}} 
                                \put(7,5){\makebox(0,0){$5$}} 
                                \put(1,1){\makebox(0,0){$6$}} 
                                \end{picture}}} 
\put(57,38){\makebox(0,0){$Q'(126453)=$}} 
\put(57,33){\makebox(0,0){\begin{picture}(8,6)(0,0) 
                                \multiput(0,0)(0,2){2}{\line(1,0){2}} 
                                \multiput(4,4)(2,0){3}{\line(0,1){2}} 
                                \multiput(0,4)(0,2){2}{\line(1,0){8}} 
                                \multiput(0,0)(2,0){2}{\line(0,1){6}} 
                                \put(1,5){\makebox(0,0){$1$}} 
                                \put(1,3){\makebox(0,0){$2$}} 
                                \put(3,5){\makebox(0,0){$3$}} 
                                \put(1,1){\makebox(0,0){$4$}} 
                                \put(5,5){\makebox(0,0){$5$}} 
                                \put(7,5){\makebox(0,0){$6$}} 
                                \end{picture}}} 
\end{picture} 
\normalsize 
\caption{Shapes $\dc_{[i;j]}$ for the permutation $\sigma=126453$ 
}
\label{fig:collat}
\end{figure} 

\begin{theorem} \label{thm:colrules} 
The shape $\dc_{[i;j+1]}$ is uniquely determined by the shapes 
$\dc_{[i-1;j]}$, $\dc_{[i;j]}$, and $\dc_{[i;j+1]}$, as follows. 
If there exists a shape $\dc \neq \dc_{[i-1;j]}$ 
such that $\dc_{[i;j]} \subsetneq \dc \subsetneq \dc_{[i-1;j+1]}$, 
then $\dc_{[i;j+1]}=\dc$. 
Otherwise $\dc_{[i;j+1]}=\dc_{[i-1;j]}$. 
\end{theorem} 

Note that in Theorem~\ref{thm:colrules}, 
a shape $\dc$ with the given properties exists if and only if 
the boxes $\dc_{[i-1;j]}-\dc_{[i;j]}$ and
$\dc_{[i-1;j+1]}-\dc_{[i-1;j]}$ 
are not adjacent to each other. 

\proof 
The second part of the theorem is clear, since the shape 
$\dc=\dc_{[i;j+1]}$ satisfies 
the condition $\dc_{[i;j]}\subsetneq \dc \subsetneq \dc_{[i-1;j+1]}$. 

Let us now assume that a shape $\dc \neq \dc_{[i-1;j]}$ 
satisfying this condition does exist. 
Suppose that, contrary to theorem's claim,
$\dc_{[i;j+1]}=\dc_{[i-1;j]}$. 
Let $A$ and $B$ be the boxes defined by 
$\{A\}=\dc_{[i-1;j]}-\dc_{[i;j]}$ 
and $\{B\}=\dc_{[i-1;j+1]}-\dc_{[i-1;j]}$. 
Denote $\tilde P=P_\sigma[i-1;j+1]$, 
$p_1=(i-1,\sigma(i-1))$ and 
$p_2=(j+1,\sigma(j+1))$. 
Then 
\[
\begin{array}{c}
\dc(\tilde P-\{p_1\})=
\dc(\tilde P-\{p_2\})=
\dc_{[i-1;j+1]}-\{B\} \,, \\[.1in]
\dc(\tilde P-\{p_1,p_2\})=\dc_{[i-1;j+1]}-\{A,B\} \,.
\end{array}
\]
Since $p_1$ and $p_2$ are minimal and maximal elements, respectively, 
of the poset~$\tilde P$, 
Theorem~\ref{thm:fominG} implies that $B$
is located either in the same column 
as~$A$, 
or in the column next to it on the right (see Figure~\ref{fig:proof2}a). 
Now let us introduce a new partial order on $\tilde P$
(denoted~$\tilde P'$)
by 
\begin{equation}
\label{eq:transversal} 
(k,\sigma(k))\leq (l,\sigma(l))
\,\Longleftrightarrow\,
\text{$k\geq l$ and $\sigma(k)\leq \sigma(l)$}\,,
\end{equation}
as in the last part of the proof of Theorem~\ref{thm:growP}. 
(So all the respective shapes get transposed.) 
Then $p_1$ and $p_2$ reverse their roles in~$\tilde P'$,
becoming maximal and minimal, respectively. 
Just as in the proof of Theorem~\ref{thm:growP}, 
Theorem~\ref{thm:fominG} implies that 
$B$ lies in the same row of $\dc_{[i-1;j+1]}$ as~$A$, 
or in the row immediately below it 
(see Figure~\ref{fig:proof2}b). 
Comparing Figures~\ref{fig:proof2}a and~\ref{fig:proof2}b, we conclude
that $B$~must be adjacent to~$A$ (see Figure~\ref{fig:proof2}c).  
This however implies that $\dc=\dc_{[i;j]} \cup \{A\}$ is 
the only shape satisfying $\dc_{[i;j]}\subsetneq \dc\subsetneq 
\dc_{[i-1;j+1]}=\dc_{[i;j]} \cup \{A,B\}$, a contradiction. 
\endproof

\begin{figure}[ht]
\centering
\setlength{\unitlength}{2.5pt}
\mbox{
   \subfigure[]{ 
        \begin{picture}(36,36)(0,0)
        \put(0, 0){\line(0,1){36}}
        \put(0,36){\line(1,0){36}}
        \put(18,18){\line(0,1){6}}
        \put(24,18){\line(1,0){6}}
        \multiput(18,12)(0,6){3}{\line(1,0){6}}
        \put(21,21){\makebox(0,0){$A$}} 
        \multiput(18,12)(1.5,0){5}{\line(0,1){6}} 
        \multiput(24,18)(1.5,0){5}{\line(0,1){18}} 
        \end{picture} 
    }
\qquad \ 
    \subfigure[]{
        \begin{picture}(36,36)(0,0)
        \put(0, 0){\line(0,1){36}}
        \put(0,36){\line(1,0){36}}
        \put(24,12){\line(0,1){12}} 
        \put(30,18){\line(0,1){6}} 
        \multiput(0 ,12)(0,6){2}{\line(1,0){24}}
        \multiput(18,18)(0,6){2}{\line(1,0){12}}
        \put(18,18){\line(0,1){6}}
        \put(21,21){\makebox(0,0){$A$}} 
        \multiput(0 ,12)(0,1.5){5}{\line(1,0){24}} 
        \multiput(24,18)(0,1.5){5}{\line(1,0){6}} 
        \end{picture}
    }
\qquad \ 
    \subfigure[]{
        \begin{picture}(36,36)(0,0)
        \put(0, 0){\line(0,1){36}}
        \put(0,36){\line(1,0){36}}
        \multiput(18,18)(6,0){2}{\line(0,1){6}}
        \multiput(18,18)(0,6){2}{\line(1,0){6}} 
        \put(21,21){\makebox(0,0){$A$}} 
        \multiput(18,12)(0,1.5){5}{\line(1,0){6}} 
        \multiput(18,12)(1.5,0){5}{\line(0,1){6}}
        \multiput(24,18)(0,1.5){5}{\line(1,0){6}} 
        \multiput(24,18)(1.5,0){5}{\line(0,1){6}}
        \end{picture}
   }
} 
\caption{Allowable locations of~$B$} 
\label{fig:proof2} 
\end{figure}
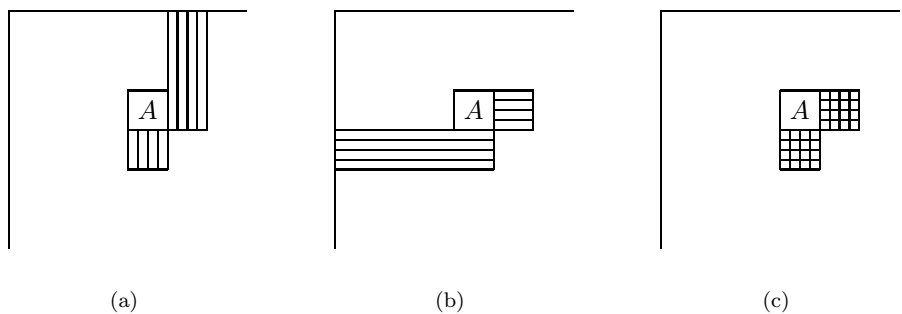

Suppose that the tableau $Q(\sigma)$ is given;
equivalently, we know the shapes 
$\dc_{[1,j]}$, for $j=1,2,\dots$;  
these are the shapes located on the upper-left side of our triangular
array. 
Recursively using the growth rules in Theorem~\ref{thm:colrules}
(moving left-to-right),
we can compute all shapes $\dc_{[i;j]}$ in the array.
In particular, we can determine 
the sequence $\dc_{[1;n]}, \dc_{[2;n]},\dots,\dc_{[n;n]}$ 
defining the tableau~$Q'(\sigma)$. 
Thus $Q'(\sigma)$ is indeed determined by $Q(\sigma)$ alone.
Since $Q'(\sigma)=Q(\sigma')$, where $\sigma'$ is $\sigma$ rotated
$180^\circ$, applying this procedure to $Q'(\sigma)$
recovers the original tableau $Q(\sigma)$. 
Thus the map $Q\mapsto Q'$ is a (shape-preserving) \emph{involution} 
on the set of standard tableaux.


Just as it was in the case of the Robinson-Schensted correspondence,
the algorithm that computes the Sch\"utzenberger involution can be ``sequentialized,'' 
and restated entirely in the language of tableaux. 
Specifically, let us first apply elementary recursion steps (based on the
rules of Theorem~\ref{thm:colrules}) to the locations adjacent to the
upper-left boundary (i.e., fix $i=2$ and take $j=2,3,\dots,n-1$, in
this order),
then to the ones adjacent to them (i.e., those with $i=3$), etc. 
Each diagonal row of shapes $\dc_{[i;i]}, \dc_{[i;i+1]},\dots,\dc_{[i;n]}$ 
is encoded by the tableau $Q_i$, which the algorithm will recursively
compute, beginning with $Q_1=Q(\sigma)$. 
The procedure that computes each tableau $Q_i$ from $Q_{i-1}$ can 
be seen to coincide with Sch\"utzenberger's ``evacuation step''~$\Delta$. 
The sequence of nested shapes of the tableaux $Q_1,\dots,Q_n$ is
recorded by the tableau $Q'(\sigma)$ (the ``evacuation tableau''). 
See Figure~\ref{fig:S}. 

Further details and references
pertaining to this fundamental combinatorial construction can be found
in \cite{fomin6, leeuwen1}; this includes applications  of
growth diagram techniques to the study of Sch\"utzenberger's \emph{jeu de
  taquin} and the proof of the \emph{Littlewood-Richardson rule.} 
(Regarding the latter, see also~\cite{FG}.) 

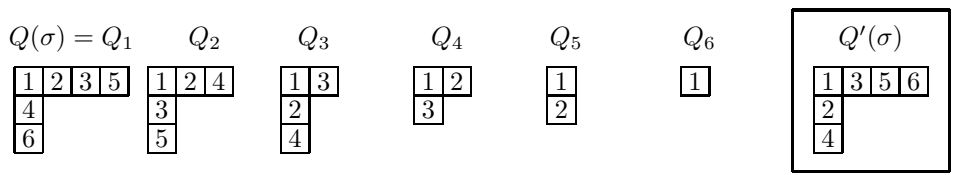
\begin{figure}[ht]
\centering 
\setlength{\unitlength}{1.8pt}
\begin{picture}(176,30)(-22,33) 
\put(149,60){\makebox(0,0){$Q'(\sigma)$}} 
\put(-19,60){\makebox(0,0){$Q(\sigma)=Q_1$}} 
\put(9  ,60){\makebox(0,0){$Q_2$}} 
\put(32 ,60){\makebox(0,0){$Q_3$}} 
\put(60 ,60){\makebox(0,0){$Q_4$}} 
\put(85 ,60){\makebox(0,0){$Q_5$}} 
\put(113,60){\makebox(0,0){$Q_6$}} 
\multiput(-28,51)(28,0){6}{\makebox(0,0){$1$}} 
\multiput(-22,51)(28,0){2}{\makebox(0,0){$2$}} 
\multiput(-16,51)(50,0){2}{\makebox(0,0){$3$}} 
\multiput(0  ,45)(56,0){2}{\makebox(0,0){$3$}} 
\multiput(28 ,45)(56,0){2}{\makebox(0,0){$2$}} 
\put(62 ,51){\makebox(0,0){$2$}} 
\put(-28,39){\makebox(0,0){$6$}} 
\put(-28,45){\makebox(0,0){$4$}} 
\put(0,  39){\makebox(0,0){$5$}} 
\put(-10,51){\makebox(0,0){$5$}} 
\put(12, 51){\makebox(0,0){$4$}} 
\put(28, 39){\makebox(0,0){$4$}} 
\put(140,51){\makebox(0,0){$1$}} 
\put(140,45){\makebox(0,0){$2$}} 
\put(146,51){\makebox(0,0){$3$}} 
\put(140,39){\makebox(0,0){$4$}} 
\put(152,51){\makebox(0,0){$5$}} 
\put(158,51){\makebox(0,0){$6$}} 
\put(149,45){\makebox(0,0){\begin{picture}(24,18)(0,0) 
                        \multiput(0 , 0)(0,6){2}{\line(1,0){6}} 
                        \multiput(0 , 0)(6,0){2}{\line(0,1){18}} 
                        \multiput(12,12)(6,0){3}{\line(0,1){6}} 
                        \multiput(0 ,12)(0,6){2}{\line(1,0){24}} 
                        \end{picture}}} 
\put(-19,45){\makebox(0,0){\begin{picture}(24,18)(0,0) 
                        \multiput(0 , 0)(0,6){2}{\line(1,0){6}} 
                        \multiput(0 , 0)(6,0){2}{\line(0,1){18}} 
                        \multiput(12,12)(6,0){3}{\line(0,1){6}} 
                        \multiput(0 ,12)(0,6){2}{\line(1,0){24}} 
                        \end{picture}}} 
\put(9,45){\makebox(0,0){\begin{picture}(24,18)(0,0) 
                        \multiput(12,12)(6,0){2}{\line(0,1){6}} 
                        \multiput(0 , 0)(0,6){2}{\line(1,0){6}} 
                        \multiput(0 , 0)(6,0){2}{\line(0,1){18}} 
                        \multiput(0 ,12)(0,6){2}{\line(1,0){18}} 
                        \end{picture}}} 
\put(37,45){\makebox(0,0){\begin{picture}(24,18)(0,0) 
                        \put(12,12){\line(0,1){6}} 
                        \multiput(0 , 0)(0,6){2}{\line(1,0){6}} 
                        \multiput(0 , 0)(6,0){2}{\line(0,1){18}} 
                        \multiput(0 ,12)(0,6){2}{\line(1,0){12}} 
                        \end{picture}}} 
\put(65,45){\makebox(0,0){\begin{picture}(24,18)(0,0) 
                        \put(0,  6){\line(1,0){6}} 
                        \put(12,12){\line(0,1){6}} 
                        \multiput(0, 6)(6,0){2}{\line(0,1){12}} 
                        \multiput(0,12)(0,6){2}{\line(1,0){12}} 
                        \end{picture}}} 
\put(93,45){\makebox(0,0){\begin{picture}(24,18)(0,0) 
                        \multiput(0 , 6)(0,6){3}{\line(1,0){6}} 
                        \multiput(0 , 6)(6,0){2}{\line(0,1){12}} 
                        \end{picture}}} 
\put(121,45){\makebox(0,0){\begin{picture}(24,18)(0,0) 
                        \multiput(0 ,12)(0,6){2}{\line(1,0){6}} 
                        \multiput(0 ,12)(6,0){2}{\line(0,1){6}} 
                        \end{picture}}} 

\thicklines
\multiput(132.5,32)(33,0){2}{\line(0,1){34}} 
\multiput(132.5,32)(0,34){2}{\line(1,0){33}} 
\end{picture} 
\caption{Sch\"utzenberger's evacuation for $\sigma=126453$.}
\label{fig:S}
\end{figure}

\end{section}

\begin{section}{Saturation and orthogonality}
\label{sec:sat}

In this section, we show that the results of Greene and
Kleitman's pioneering papers~\cite{GK1, greene3} 
on what they called ``saturated'' families of chains or antichains 
can be viewed as simple corollaries of one master theorem, the
Duality Theorem for Finite Posets (Theorem~\ref{thm:dual}). 
(It should be noted that historically, 
the sequence of events was different: 
the main saturation result in \cite{GK1} was a principal tool in
Greene's original proof~\cite{greene3} of the Duality Theorem.)

In what follows, ``chain family'' always means a collection of
disjoint chains. 
(The term ``chain $k$-family'' will emphasize that there are $k$
chains in this collection.)
If, in addition, these chains cover the whole poset,
they are said to form its ``chain partition.''
A chain $k$-family is \emph{maximal} if it covers the maximal possible
number of elements of~$P$. 
The same conventions will apply to antichain families and partitions.

Let  $\mc$ be a chain partition of a finite poset~$P$. 
Since a chain may intersect an antichain 
in at most one element, the total size of any $k$ disjoint antichains
$A_1,\dots,A_k\subset P$ is
bounded from above by a quantity that only depends on~$\mc$:   
\begin{equation} 
\label{eq:saturated}
\sum_{i=1}^k |A_i|=\sum_{C\in \mc}\sum_{i=1}^k |A_i\cap C|
\leq \sum_{C\in \mc} \min \{|C|,k\} \,. 
\end{equation} 
Similarly, the total size of any family of $k$ disjoint chains is at most 
$\displaystyle\sum_{A\in \ma}\min \{ |A|,k\}$, given an antichain
partition $\ma$ of~$P$. 
A chain (resp.\ antichain) partition $\mc$ (resp.~$\ma$) 
is said to be \emph{$k$-saturated} 
if the upper bound described above is achieved for some 
disjoint family of $k$ antichains (resp.\ chains).  



\begin{theorem} \label{thm:kk+1chain} \cite{GK1, greene3} 
For each $k$, there exists a chain (resp.\ antichain) partition  of
$P$ which is simultaneously $k$-saturated and $(k+1)$-saturated.
\end{theorem}

The chain version of Theorem~\ref{thm:kk+1chain} is due to 
Greene and Kleitman~\cite{GK1};
another proof was later given by H.~Perfect \cite{perfect2},
using the idea of M.~Saks~\cite{saks1}.
The antichain counterpart was obtained by Greene~\cite{greene3}. 

It is already quite non-trivial to show 
that a $k$-saturated partition exists for every~$k$. 
(See a short proof in~\cite{saks1}.) 
For $k=1$, the existence of a $1$-saturated chain
partition is equivalent to Dilworth's theorem.



For many classes of posets
(e.g., Boolean algebras), 
there always exists a chain partition which is $k$-saturated for
\emph{all~$k$}. 
This is however false in general. 
The poset in Figure~\ref{fig:dual} provides a counterexample:
the only $1$-saturated chain partition is $\{bdf, ace\}$,
which is not $3$-saturated, as the right-hand side of 
(\ref{eq:saturated}) (with $k=3$) is~$n=6$ in this case, 
while the left-hand side is at most~5 for a family of 3~antichains. 

The derivation of Theorem~\ref{thm:kk+1chain} 
from Theorem~\ref{thm:dual} given below 
employs the concept of \emph{orthogonality}, which 
plays a major role in the proofs of the Duality Theorem given 
in~\cite{fomin1, frank},
as well as in directed graph generalizations developed by
S.~Felsner~\cite{felsner}. 

\begin{definition}
\label{def:ortho}
In a finite poset~$P$, 
a chain family $\mc=\{C_1,\dots,C_l\}$ and an antichain 
family $\ma=\{A_1,\dots,A_k\}$ are called \emph{orthogonal} if 
\begin{align}
\label{eq:orth1} 
&P=C_1 \cup \cdots\cup C_l \cup A_1\cup\cdots\cup A_k\,;\\
\label{eq:orth2} 
&\text{$C_i\cap A_j\neq \phi$ for all $i=1,\dots,l$, $j=1,\dots,k$.}
\end{align}
\end{definition}


The notion of orthogonality can be reformulated as follows. 

\begin{lemma}
\label{lem:A+B=n+kl}
In a finite poset~$P$, of cardinality~$n$,
a chain family $\mc=\{C_1,\dots,C_l\}$ and 
an antichain family $\ma=\{A_1,\dots,A_k\}$ are orthogonal if and only if 
\[
\sum_i|C_i| + \sum_j|A_j| = n+kl\,.
\]
If the families $\mc$ and $\ma$ are orthogonal, then they are both
maximal. 
\end{lemma}

\proof
Since a chain and an antichain intersect in at most one element, 
we have
\begin{equation*} 
\label{eq:incl-excl}
\sum_i|C_i| + \sum_j|A_j| 
= \left|\bigcup C_i\right| + \left|\bigcup A_j\right| 
=\left|\bigcup C_i \cup \bigcup A_j\right|
  +\sum_{i,j}|C_i\cap A_j|
\leq n+kl \,,
\end{equation*}
with equality holding if and only if the conditions 
(\ref{eq:orth1})--(\ref{eq:orth2}) are satisfied. 
This proves the first part of the lemma. 
The second part follows as well, for if $\mc$ and $\ma$ were
orthogonal but not maximal, then the last inequality would be violated
by the maximal families of the corresponding sizes. 
\endproof


We next explain the connection between saturation and orthogonality. 

\begin{lemma}
\label{lem:ortho-sat}
Let $\mc=\{C_1,\dots,C_l\}$ and $\ma=\{A_1,\dots,A_k\}$ 
be orthogonal chain and antichain families, respectively. 
Let $\mc^+$ be the chain partition consisting of $\mc$ 
together with the elements of
the complement $P-\cup C_i$ viewed as single-element chains. 
Then $\mc^+$ is a $k$-saturated chain partition. 
The analogously defined antichain partition $\ma^+$ is $l$-saturated. 
\end{lemma}

\proof
Since each chain $C_i$ intersects all $k$ (disjoint)
antichains in $\ma$, we conclude that $|C_i|\geq k$. 
Then the right-hand side of (\ref{eq:saturated}),
with $\mc$ replaced by $\mc^+$, 
is equal to $kl+n-\sum|C_i|$. 
By Lemma~\ref{lem:A+B=n+kl},
the latter expression equals $\sum |A_j|$, as desired. 
\endproof

We are now prepared to prove Theorem~\ref{thm:kk+1chain}
(using Theorem~\ref{thm:dual}). 

\subseqnoref{Proof of Theorem~\ref{thm:kk+1chain}} 
We will prove the chain partition version, as the proof of the
antichain counterpart is completely analogous. 
Denote $\dc=\dc(P)$, and let $l$ be 
uniquely defined by the condition that $\dc_l>k$ 
while $\dc_{l+1}\leq k$; 
in other words,  the points $(k,l)$ and $(k+1,l)$ lie on the outer
boundary of the shape $\dc$  
(see Figure~\ref{fig:orth}). 
Let $\ma$ be a maximal antichain $k$-family, 
and $\mc$ a maximal chain $l$-family. 
Thus $\mc$ covers $\dc_1+\cdots+\dc_l$ elements,
while $\ma$ covers $\dc'_1+\cdots+\dc'_k$ elements,
where $\dc'_1,\dc'_2,\dots$ are the column lengths of~$\dc$. 
Then 
$\sum_i|C_i| + \sum_j|A_j| = \dc_1+\cdots+\dc_l + \dc'_1+\cdots+\dc'_k
= n+kl$;
so by Lemma~\ref{lem:A+B=n+kl}, 
the families $\mc$ and $\ma$ are orthogonal. 
Hence by Lemma~\ref{lem:ortho-sat}, 
the chain partition $\mc^+$ obtained by adding singletons to $\mc$ is
$k$-saturated. By the same token, it is $(k+1)$-saturated, and the
theorem is proved. 
\endproof

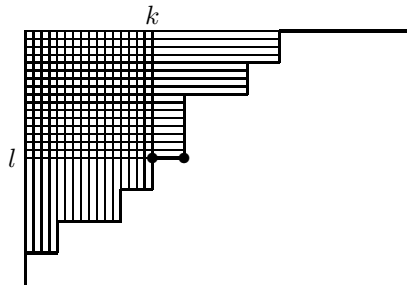
\begin{figure}[ht]
\centering
\setlength{\unitlength}{2pt}
\begin{picture}(75,45)(-3,-3)
\multiput(0,0)(1.5,0){5}{\line(0,1){42}}
\multiput(6,6)(1.5,0){9}{\line(0,1){36}}
\multiput(18,12)(1.5,0){5}{\line(0,1){30}}
\multiput(0,18)(0,1.5){8}{\line(1,0){30}}
\multiput(0,30)(0,1.5){4}{\line(1,0){42}}
\multiput(0,36)(0,1.5){5}{\line(1,0){48}}
\put(24,18){\circle*{2}}  
\put(30,18){\circle*{2}} 
\put(24,45){\makebox(0,0){$k$}}
\put(-2.5,18){\makebox(0,0){$l$}}
\thicklines
\put(0,-6){\line(0,1){6}}
\put(0,0){\line(1,0){6}}
\put(6,0){\line(0,1){6}} 
\put(6,6){\line(1,0){12}}
\put(18,6){\line(0,1){6}} 
\put(18,12){\line(1,0){6}}
\put(24,12){\line(0,1){6}}
\put(24,18){\line(1,0){6}}
\put(30,18){\line(0,1){12}}
\put(30,30){\line(1,0){12}} 
\put(42,30){\line(0,1){6}} 
\put(42,36){\line(1,0){6}} 
\put(48,36){\line(0,1){6}}
\put(48,42){\line(1,0){24}}
\end{picture} 
\caption{Orthogonality}
\label{fig:orth}
\end{figure}

\end{section}

\section{Nilpotent matrices}
\label{sec:nilpotent}

We will now discuss an important interpretation of the shape $\dc(P)$ 
in terms of sizes of Jordan blocks for a ``typical'' element of the
incidence algebra of the poset~$P$ 
(see Theorem~\ref{thm:Jordan-generic} below). 
This connection was discovered independently by Michael
Saks~\cite{saks2, saks3} and Emden Gansner~\cite{gansner}, 
and actually extends to the broader setting of acyclic directed graphs. 
In this section, we follow the general plan of
\cite[pp.~429--431]{gansner}, 
restricted to the case of posets.  

Let us fix a labelling identifying the poset $P$ with the set
$\{1,\dots,n\}$. 
It will be convenient to assume that our labelling is a linear
  extension, i.e., larger elements receive larger labels. 
The \emph{incidence algebra} $I(P)$ can be defined as 
the set of complex matrices $M$ such that 
$M_{ij}\neq 0$ implies $i\leq j$ in~$P$
(cf. \cite[3.6]{stanley1}). 
In particular, all these matrices are upper-triangular.  
A nilpotent element $M\in I(P)$ (i.e., such that $M_{ii}=0$ for
all~$i$) is called \emph{generic} if the entries 
$M_{ij}\,$, $i\stackrel{P}<j$ are independent transcendentals
(over~$\mathbb{Q}$); in particular, all of these entries must be
nonzero. 
Figure~\ref{fig:generic}a shows a generic nilpotent element of $I(P)$
for the poset and the labelling shown in Figure~\ref{fig:PQ}a. 
(The values $M_{ij}$ are presumed algebraically independent.) 

The Jordan canonical form of an $n\times n$ complex matrix $M$
consists of a number 
of Jordan blocks of sizes $n_1\geq n_2\geq \cdots $, 
the \emph{invariants} of~$M$. 
The partition $J(M)=(n_1,n_2,\dots )$ of $n$ is called 
the \emph{Jordan partition} of~$M$.
We append zeroes at the end of the sequence $n_1,n_2,\dots$, so the number
$n_i$ makes sense for any~$i>0$. 
 
\begin{theorem} {\rm \cite{gansner,saks2}}
\label{thm:Jordan-generic}
The Jordan partition of a generic nilpotent element of the incidence
algebra $I(P)$ is $\dc(P).$ 
\end{theorem}

To illustrate Theorem~\ref{thm:Jordan-generic}, 
Figure~\ref{fig:generic}b shows the Jordan canonical
form of the matrix on the left. 
Thus the Jordan partition in this case is $(4,2)$, 
in agreement with Figure~\ref{fig:dual}.

\begin{figure}[ht]
\centering
\setlength{\unitlength}{3pt}
\mbox{
\subfigure[A generic nilpotent element $M\in I(P)$]{
        $\!\!\!\left[\begin{array}{cccccc}
              \quad\!\!\!0\quad\!\!\! & M_{12} & M_{13} & \quad\!\!\!0\quad\!\!\! & M_{15} & M_{16} \\[.1in]
              0 & 0 & M_{23} & 0      & M_{25} & M_{26} \\[.1in]
              0 & 0 & 0      & 0      & 0      & 0      \\[.1in]
              0 & 0 & 0      & 0      & M_{45} & M_{46} \\[.1in]
              0 & 0 & 0      & 0      & 0      & M_{56} \\[.1in]
              0 & 0 & 0      & 0      & 0      & 0      
          \end{array}\right]$
}

    \subfigure[The Jordan canonical form of $M$]{
        $\left[\begin{array}{cccc||cc}
              \quad\!\!\!0\quad\!\!\! & \quad\!\!\!1\quad\!\!\! & \quad\!\!\!0\quad\!\!\! & \quad\!\!\!0\quad\!\!\! & \quad\!\!\!0\quad\!\!\! & \quad\!\!\!0\quad\!\!\! \\[.1in]
              0 & 0 & 1 & 0 & 0 & 0 \\[.1in]
              0 & 0 & 0 & 1 & 0 & 0 \\[.1in]
              0 & 0 & 0 & 0 & 0 & 0 \\[.07in]
\hline
\hline
              0 & 0 & 0 & 0 & 0 & 1 \\[.1in]
              0 & 0 & 0 & 0 & 0 & 0 
          \end{array}\right]$        } 
}
\caption{Theorem~\ref{thm:Jordan-generic}}
\label{fig:generic}
\end{figure}

A more general statement, which describes the Jordan partition of
a generic nilpotent matrix with a fixed pattern of zeroes, 
can be obtained from a result by S.~Poljak~\cite{poljak}. 

\subseqnoref{Proof of Theorem~\ref{thm:Jordan-generic}}
Let $M$ be a generic nilpotent element in~$I(P)$, 
and let $n_1\geq n_2\geq \cdots$ be its invariants. 
Let $x$ be a formal variable,
and let $p_k(xI-M)$ denote the greatest common divisor, with leading
coefficient 1, of all $k\times k$ minors of the matrix $xI-M$. 
For $k\leq 0$, we set $p_k(xI-M)=1$. 
Note that $p_n(xI-M)=x^n$. 
An example is given below:
\[
M=\left[\begin{array}{ccc}
              0 & M_{12} & 0 \\[.1in]
              0 & 0 & 0 \\[.1in]
              0 & 0 & 0 
          \end{array}\right]
\qquad
xI-M=\left[\begin{array}{ccc}
              x & -M_{12} & 0 \\[.1in]
              0 & x & 0 \\[.1in]
              0 & 0 & x 
          \end{array}\right]
\qquad
\begin{array}{l}
p_1(xI-M)=1\\[.1in]
p_2(xI-M)=x\\[.1in]
p_3(xI-M)=x^3
\end{array} 
\]
We will need the following basic linear-algebraic result
(see, e.g., \cite[6.43]{shilov}), which does not require the
assumption of genericity. 

\begin{lemma} \label{lem:gansner22} 
For $k\leq n$, we have
$p_k(xI-M)=x^{d_k}$,
where $d_k=n-\sum_{i\leq n-k} n_i\,$. 
\end{lemma}


If $a_1(P)=n$, then $P$ is an antichain, $M=0$, and $n_k=\dc_k=1$ for
all $k\geq 1$. 
Assume therefore that $a_1(P)<n$. 
Let $p_k(xI-M)=x^{d_k}$, as in Lemma~\ref{lem:gansner22}.
In order to prove that partitions $(n_1,n_2,\dots)$ and $\dc(P)$
coincide, we need to show that $d_{n-k}=n-c_k(P)$ for all $k\geq 0$. 
This holds trivially for $k=0$, 
so assume $k>0$. 
It follows from the definition that 
$p_{n-k}(xI-M)=x^{d_{n-k}}$ is the smallest power 
of $x$ appearing as a term in any $(n-k)$-minor of $xI-M$. 
Since $M$ is generic, $d_{n-k}$ is the smallest number of diagonal
entries appearing in any collection of $n-k$ nonzero entries of $xI-M$, 
no two in the same row or column. 

The claim $d_{n-k}=n-c_k(P)$ now becomes a purely combinatorial
statement, which we will now verify. 
(This statement can be reformulated and proved using the 
network construction due to A.~Frank, 
to be introduced and studied in Section~\ref{sec:frank's-proof}. 
To keep this part of our presentation self-contained, 
an independent proof follows.) 

Assume that $k\leq a_1(P)$, and let $\mc=\{C_1,...,C_k\}$ be a maximal
chain $k$-family in~$P$. 
For each chain $C_i=(p_1<p_2<\cdots<p_l)$ of $\mc$ containing $l\geq 2$ 
elements, consider the $l-1$ entries of the matrix 
$xI-M$ located in positions
$(p_1,p_2)$, $(p_2,p_3)$, \dots, $(p_{l-1},p_l)$. 
The total number of such entries is $c_k(P)-k$. 
Since the chains of $\mc$ are disjoint, 
no two entries occupy the same row or column. 
Throw in the $n-c_k(P)$ diagonal entries $(xI-M)_{p,p}$ corresponding 
to the elements $p\in P$ not covered by~$\mc$. 
In total, all these entries number $n-k$, of which $n-c_k(P)$ are
equal to $x$. 
Hence $d_{n-k}\leq n-c_k(P)$. 

Conversely, consider $n-k$ entries  
no two of which are in the same row or column, 
and suppose that $d_{n-k}$ of them are diagonal entries (each equal
to~$x$). 
The remaining $n-k-d_{n-k}$ entries correspond to a disjoint collection 
$\{C_1,...,C_l\}$ of chains, each containing at least two elements.  
Together these chains cover $n-k-d_{n-k}+l$ elements. 
With the elements corresponding to the diagonal entries, 
they total $n-k+l\leq n$ elements. 
Hence $l\leq k$. 
There are $k-l+d_{n-k}\geq 0$ elements not covered by 
the chains $C_1,...,C_l$. 
Choose $k-l$ of these, say $p_1,...,p_{k-l}$ and form the chain $k$-family 
$\{C_1,...,C_l,\{p_1\},...,\{p_{k-l}\} \}$. 
Then $n-c_k(P)\leq n-(|\cup C_i|+k-l)=d_{n-k}$. 
We thus proved that $d_{n-k}=n-c_k(P)$ for all $k\leq a_1(P)$. 
Hence $n_i=\dc_i(P)$ for all $i\leq a_1(P)$, and therefore for all
$i\geq 0$. 
\endproof

It is well known (and easy to prove)
that the closure of the set of nilpotent matrices
with Jordan partition $\dc$ consists of all nilpotent matrices whose 
Jordan partition $\mu$ is $\leq\dc$ with respect to the
\emph{dominance} order (i.e., 
$\mu_1+\cdots+\mu_i\leq \dc_1+\dots+\dc_i$ for all~$i$).
It then follows from Theorem~\ref{thm:Jordan-generic}  
that $\dc(P)$ dominates the Jordan partition of \emph{any}
nilpotent element in~$I(P)$. 
A direct proof of this statement was given in~\cite{saks2}.

In the case of permutation posets, Theorem~\ref{thm:Jordan-generic} 
leads to an important geometric interpretation of the
Robinson-Schensted correspondence discovered by
Robert Steinberg~\cite{steinberg}.
We will now briefly (and informally) describe the main combinatorial
ingredients of Steinberg's construction, 
trying to keep our presentation elementary;
see the original paper~\cite{steinberg}
or Marc van Leeuwen's insightful exposition~\cite{leeuwen2} for
further details, 
and the work of van Leeuwen~\cite{mvl-thesis} and Itaru Terada 
\cite{terada} for generalizations of Steinberg's construction. 

Let $e_1,\dots,e_n$ be the standard linear basis in~$\mathbb{C}^n$. 
Let $\sigma$ be a permutation of~$[n]$, and let $E=(E_1,\dots,E_n)$
and $F=(F_1,\dots,F_n)$ be the flags of subspaces defined by
\begin{eqnarray}
\label{eq:EF}
\begin{array}{ll}
E_1={\rm span}(e_1), & F_1={\rm span}(e_{\sigma(1)}), \\[.1in]
E_2={\rm span}(e_1, e_2), & F_2={\rm span}(e_{\sigma(1)}, e_{\sigma(2)}), \\[.1in]
\cdots & \cdots \\[.1in]
E_n={\rm span}(e_1,\dots,e_n), \qquad
& F_n={\rm span}(e_{\sigma(1)},\dots,e_{\sigma(n)}). 
\end{array}
\end{eqnarray}
(Thus $E$ and $F$ are \emph{in position~$\sigma$} with respect to each
other.) 
It is straightforward to verify that the incidence algebra $I(P_\sigma)$ 
of the permutation poset $P_\sigma$ is exactly the set of matrices
which fix each of the subspaces $E_i$ and $F_j$ (i.e., $I(P_\sigma)$
is the common stabilizer of $E$ and~$F$). 
We thus obtain the following corollary of Theorem~\ref{thm:Jordan-generic}.

\begin{corollary}\cite{steinberg}
The Robinson-Schensted shape $\dc(P_\sigma)=|P(\sigma)|=|Q(\sigma)|$
of a permutation $\sigma$ can be defined as the Jordan partition of a
generic nilpotent matrix that fixes two flags in relative
position~$\sigma$ with respect to each other. 
\end{corollary}

Viewing a generic nilpotent element $M\in I(P_\sigma)$ 
as a matrix of a linear transformation, consider the restriction of
this transformation to an invariant subspace~$E_i\,$.
The matrix $M|_{E_i}$ of this restriction is the principal submatrix
of $M$ obtained by 
taking  the first $i$ rows and the first $i$ columns.
This submatrix is obviously a generic nilpotent element of the
incidence algebra  $I(P_\sigma(n,i))$, where $P_\sigma(n,i)$ is
the subposet of $P_\sigma$ formed by the elements with labels
$1,\dots,i$ (cf.\ (\ref{eq:P(i,j)})). 
Comparing this observation to the description of the ``insertion
tableau'' $P(\sigma)$ given in Section~\ref{sec:schensted}, we 
arrive at the following conclusion. 

\begin{corollary}
\label{cor:steinberg2} 
\cite{steinberg}
The Robinson-Schensted correspondence
$\sigma\mapsto(P(\sigma),Q(\sigma))$ has the following geometric
interpretation. 
The tableau $P(\sigma)$ records the growth of Jordan partitions for
the restrictions $M|_{E_i}\,$,
where $M$ is a generic nilpotent transformation that fixes two flags
$E$ and $F$ in relative position~$\sigma$ with respect to each other
(cf. {\rm(\ref{eq:EF})}). 
Analogously, the tableau $Q(\sigma)$ records  the growth of Jordan
partitions for the restrictions~$M|_{F_i}\,$. 
\end{corollary}

For example, the Jordan partitions of the principal submatrices of the
matrix in Figure~\ref{fig:generic}a form the tableau in
Figure~\ref{fig:PQ}a (cf.\ the right edge of Figure~\ref{fig:growP}).  

Let $Fl_M$ denote the variety of flags fixed by a given nilpotent
matrix $M$ with Jordan partition~$\lambda$.
Steinberg shows that the irreducible components of $Fl_M$ are all of
the same dimension, and are labelled
by the standard tableaux $T$ of shape~$\lambda$ in the way suggested by
Corollary~\ref{cor:steinberg2}, namely:
the irreducible component labelled by~$T$ is the closure
of the set of flags $L=(L_1,\dots,L_n)\in Fl_M$ such that the Jordan
types of the restrictions $M|_{L_i}$ are recorded by~$T$.
Further exploring this approach, M.~van Leeuwen \cite{leeuwen1,leeuwen2}
gave detailed geometric
interpretations of the tableau-theoretic constructions of
Sections~\ref{sec:schensted}--\ref{sec:schutzenberger}. 




\section{Network flow preliminaries}
\label{sec:flow-prelim}

We will now recall the minimum cost flow algorithm of Ford and
Fulkerson~\cite[III.3]{FF}. 
Further details 
can be found in~\cite{FF}, 
as well as in numerous textbooks. 

Suppose we are given a \emph{network}
with the underlying directed graph $G=(V,E)$, the \emph{source}~$s$, 
the \emph{sink}~$t$, the \emph{cost} function $a:E\to\{0,1,\dots\}$,
and the \emph{capacity} function $c:E\to\{1,2,\dots\}$. 
A \emph{flow} 
in this network is a function $f:E\to\{0,1,\dots\}$ 
satisfying the \emph{conservation law} 
$
\displaystyle\sum_{e=(x,\cdot)\in E} f(e)=\sum_{e=(\cdot,x)\in E} f(e)
$
(for $x\notin\{s,t\}$) and the \emph{capacity
  restrictions} $0\leq f(e)\leq c(e)$. 
Ford and Fulkerson's algorithm (see Algorithm~\ref{alg:MC} below)
finds a flow that has the given total \emph{value} 
\[
v=\sum_{e=(s,\cdot)\in E} f(e)=\sum_{e=(\cdot,t)\in E} f(e)
\]
and minimizes the total cost
$\sum_{e\in E} a(e)f(e)$.  

A \emph{potential} is a function $\pi:V\to\{0,1,\dots\}$
satisfying the condition $\pi (s)=0$. 
The number $p=\pi (t)$ is the \emph{value} of~$\pi$.

\begin{theorem} {\rm \cite[p.~115]{FF}}
\label{th:FFduality}
Assume that the flow $f$ and the potential $\pi$ satisfy
the following conditions:
\begin{eqnarray}
\label{eqn:mc1}
& 0\leq \pi (x)\leq p, 
& \text{for any vertex $x\in V$}\,;\\
\label{eqn:mc2} 
& \pi (y)-\pi (x)<a(x,y) \ \ \Rightarrow \ \ f(x,y)=0,\ \qquad
& \text{for any edge $(x,y)\in E$}\,;\\ 
\label{eqn:mc3} 
& \pi (y)-\pi (x)>a(x,y) \ \ \Rightarrow \ \ f(x,y)=c(x,y),
& \text{for any edge $(x,y)\in E$}. 
\end{eqnarray}
Then $f$ has minimal cost among all flows of the same value. 
\end{theorem}

\begin{algorithm}[The minimal cost flow algorithm~\cite{FF}] 
\label{alg:MC} {\rm 
\ \par 
\medskip
\noindent
\emph{Initial data:} 
flow $f$ and potential $\pi$ satisfying  
{\rm (\ref{eqn:mc1})--(\ref{eqn:mc3})}. 
In particular, we may set $\pi(x)=0$ for all $x\in V$,  
and $f(e)=0$ for all $e\in E$. 

\subseqnoref{MC1} 
Let $G'=(V,E')$ be the directed graph on the same set of vertices~$V$,
with 
\begin{align}
\notag 
E' =\ \ & \{(x,y):(x,y)\in E,\ \pi (y)-\pi (x)=a(x,y),\
f(x,y)<c(x,y)\}\\
\notag
\cup &\{(y,x):(x,y)\in E,\ \pi (y)-\pi (x)=a(x,y),\ f(x,y)>0\}\,. 
\end{align}
Let $X\subseteq V$ be the set of vertices $x$
for which a path from $s$ to $x$ exists in~$G'$. 
If $t\in X$, then go to {\bf MC2a}. 
Otherwise go to {\bf MC2b}.

\subseqnoref{MC2a} 
Let $\mathcal M$ be a path in $G'$ from $s$ to~$t$. 
Increase the flow $f$ along $\mathcal M$ by~1. 
Proceed to {\bf MC3}.

\subseqnoref{MC2b} 
Increase the potential $\pi (x)$ of each vertex $x\in V-X$ by~1.

\subseqnoref{MC3} 
If the flow is not maximal, return to {\bf MC1}; otherwise stop. 
(The maximality of the flow can be detected using similar techniques;
we will not discuss them here.) 
} 
\end{algorithm} 

We will identify each stage of Algorithm~\ref{alg:MC} 
(more precisely, a moment after/before executing {\bf MC2a}/{\bf MC2b}) 
by the corresponding pair of values~$(p,v)$. 

\begin{theorem} {\rm \cite{FF}}
\label{th:FF}
Algorithm~\ref{alg:MC} terminates by arriving at a maximal flow. 
Conditions {\rm (\ref{eqn:mc1})-(\ref{eqn:mc3})} are preserved 
throughout; thus the current flow $f$ has
minimal cost among all flows of the same value. 
If the algorithm starts with zero initial data,
it produces, at different stages, the minimal cost flows of
all possible values~$v$. 

\nopagebreak[4]
Each time the flow increases by~1, its cost 
increases by the current potential value~$p$. 
\end{theorem}


\section{Frank's network. Proof of Theorem~\ref{thm:dual}} 
\label{sec:frank's-proof}

\begin{figure}[ht]
\centering
\setlength{\unitlength}{2.6pt}
\begin{picture}(50,71)(0,0)
  \multiput(0,20)(10,0){6}{\line(0,1){30}}

  \multiput(0,20)(10,0){6}{\circle*{1.5}}
  \multiput(0,50)(10,0){6}{\circle*{1.5}}
  \put(0,70){\circle*{1.5}}
  \put(50,0){\circle*{1.5}}

  \put(0,70){\line(0,-1){20}}
  \put(0,70){\line(1,-2){10}}
  \put(0,70){\line(1,-1){20}}
  \put(0,70){\line(3,-2){30}}
  \put(0,70){\line(2,-1){40}}
  \put(0,70){\line(5,-2){50}}

  \put(50,0){\line(0,1){20}}
  \put(50,0){\line(-1,2){10}}
  \put(50,0){\line(-1,1){20}}
  \put(50,0){\line(-3,2){30}}
  \put(50,0){\line(-2,1){40}}
  \put(50,0){\line(-5,2){50}}

  \put(0,50){\line(1,-3){10}}
  \put(10,50){\line(1,-3){10}}
  \put(30,50){\line(1,-3){10}}
  \put(40,50){\line(1,-3){10}}

  \put(0,50){\line(2,-3){20}}
  \put(30,50){\line(2,-3){20}}

  \put(10,50){\line(1,-1){30}}

  \put(0,50){\line(4,-3){40}}
  \put(10,50){\line(4,-3){40}}

  \put(0,50){\line(5,-3){50}}

  \put(-3,70){\makebox(0,0){$s$}}
  \put(53,0){\makebox(0,0){$t$}}

  \put(-3,49){\makebox(0,0){$x_f$}}
  \put( 7,49){\makebox(0,0){$x_d$}}
  \put(17,49){\makebox(0,0){$x_b$}}
  \put(27,49){\makebox(0,0){$x_e$}}
  \put(37,49){\makebox(0,0){$x_c$}}
  \put(47,49){\makebox(0,0){$x_a$}}

  \put(-3,19){\makebox(0,0){$y_f$}}
  \put( 7,19){\makebox(0,0){$y_d$}}
  \put(17,19){\makebox(0,0){$y_b$}}
  \put(27,19){\makebox(0,0){$y_e$}}
  \put(37,19){\makebox(0,0){$y_c$}}
  \put(47,19){\makebox(0,0){$y_a$}}

\end{picture}
\caption{Frank's network}
\label{fig:frank-network}
\end{figure}
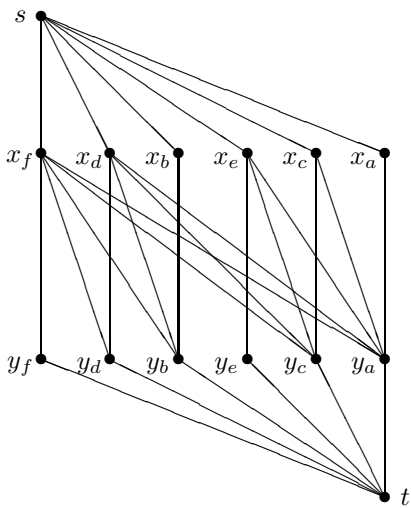

The proof of the Duality Theorem presented below is due to
A.~Frank~\cite{frank} (reproduced in~\cite{engel}). 
The main tool is an application of Algorithm~\ref{alg:MC}
to a certain network associated with the poset~$P$. 
The underlying graph $G=(V,E)$ of this network is obtained by adjoining
a source $s$ and a sink $t$ to \emph{two copies} of~$P$,
and introducing the edges, as follows: 
\begin{align*}
& V=\{ s,t\} 
\cup \{x_p: p\in P\} 
\cup \{y_p: p\in P\}\,, \\
& E=\{ (s,x_p): p\in P\} 
\cup \{(x_p,y_{p'}):p\geq p'\text{ for }p,p'\in P\} 
\cup \{(y_p,t): p\in P\}\,. 
\end{align*}
All edge capacities $c(e)$ are equal to~1,  
and the cost function is defined by 
\[
\notag a(e) = 
        \begin{cases}
        1 & \text{if $e=(x_p,y_p)$, $p\in P$;} \\
        0 & \text{otherwise.}
        \end{cases}
\]
See Figure~\ref{fig:frank-network}, which shows this network 
for the poset in Figure~\ref{fig:dual}.

The first key observation of Frank's was that the flows and potentials in
this network correspond naturally to chain and antichain families in
the underlying poset~$P$. 
Let us explain. 

With any integral flow $f$ in $G$, let us associate a chain family
$\mc(f)$ constructed as follows. 
Let the flow cost and value be equal to~$d$ and $v$, respectively. 
Thus the flow saturates $v$ edges of the form $(x_p,y_{p'})$,   
with no two of them sharing vertices, due to the capacity
restrictions. 
Among these $v$ edges, there are $d$ edges with $p=p'$;  
the remaining $v-d$ edges $(x_p,y_{p'})$, $p>p'$, 
``project'' onto two-element chains $\{p>p'\}$ in~$P$. 
These chains glue naturally into a chain family $\mc(f)$,
which covers $n-d$ elements of~$P$, and  
consists of $(n-d)-(v-d)=n-v$ chains
(here, as before, $n$ denotes the cardinality of~$P$). 

We also associate an antichain family $\ma(\pi)$ 
with an arbitrary potential function $\pi$ in the network described above. 
Let $p$ be the value of~$\pi$. 
For $i=1,\dots,p$, we let $P_i=\{ p:\pi (x_p)<\pi (y_p)=i\}$, 
and define $A_i$ to be the set of the maximal elements of~$P_i$. 
We then set $\ma(\pi)=\{A_1,\dots,A_p\}$. 


Let us apply Algorithm~\ref{alg:MC} (with zero initial data) 
to the network associated to~$P$.  
The step \textbf{MC2a} involves an arbitrary choice of a
breakthrough path~$\mathcal{M}$, so let us fix a particular sequence
of such choices once and for all. 
Consider the flow $f$ and the potential $\pi$ obtained at 
the stage $(p,v)$ of the algorithm,  
and let $\mc_{n-v}=\mc(f)$ and $\ma_p=\ma(\pi)$ 
be the corresponding chain and antichain families.  
This notation is indeed unambiguous, since $\mc_{n-v}$ 
depends only on $v$, while $\ma_p$ depends only on~$p$
(assuming that the sequence of path choices has been fixed);
this is because each execution of the loop modifies the flow or the
potential, but not both, and furthermore both $p$ and $v$ weakly
increase in the process. 

For the poset $P$ in Figure~\ref{fig:dual}, 
this procedure consecutively generates: 
\[
\begin{array}{ll}
\mc_2= \{ace,bdf \}, \qquad \\
& \ma_1 = \{ ab\}, \\
& \ma_2 = \{ ab,de\}, \\
\mc_1= \{acdf \},    \qquad \\
& \ma_3 = \{ a,bc,de\}, \\
& \ma_4 = \{ a,bc,d,ef\} . 
\end{array}
\]

The following crucial lemma is due to A.~Frank. 

\begin{lemma} \label{lem:frank} \cite{frank} 
At any stage $(p,v)$ of the minimal cost flow algorithm,
the families $\mc_{n-v}$ and $\ma_p$ are orthogonal.
\end{lemma}

\proof 
For Frank's network, conditions (\ref{eqn:mc1})--(\ref{eqn:mc3}) 
are restated as follows: 
\begin{eqnarray} 
\label{eqn:defcri0} & 0\leq \pi (x)\leq p\,, 
& \text{for any vertex $x$}\,; \\
\label{eqn:defcri1} & f(x,y)=1 \Rightarrow \pi (y)-\pi (x)\geq a(x,y)\,,
& \text{for any edge $(x,y)$}\,; \\
\label{eqn:defcri2} & f(x,y)=0 \Rightarrow \pi (y)-\pi (x)\leq a(x,y) \,,
& \text{for any edge $(x,y)$}.
\end{eqnarray}
Let $p\in P-\cup \mc_{n-v}$, 
where $\cup \mc_{n-v}$ denotes the union of all chains
in~$\mc_{n-v}\,$. 
By construction, $f(x_p,y_p)=1$.
Condition (\ref{eqn:defcri1}) implies that $\pi (y_p)\geq 1+\pi (x_p)$;  
hence $p\in P_i\,$, where $i=\pi (y_p)$. 
(Here we retain the notation introduced in the definition of~$\ma(\pi)$.) 
Suppose for a moment that $p\notin A_i\,$, 
i.e., $p<p'$ for some $p'\in P_i\,$. 
Then $\pi (y_p)=\pi (y_{p'})=i>\pi (x_{p'})$, 
by the definition of~$P_i\,$. 
Now condition (\ref{eqn:defcri2}) implies $f(x_{p'},y_p)=1$,
which contradicts the capacity restrictions. 
Thus our assumption was false, that is, $p\in A_i\,$, 
and we have proved that $(\cup \mc_{n-v})\cup (\cup \ma_p)=P$
(cf.~(\ref{eq:orth1})).

Let $C=\{ p_1>\cdots>p_b\}$ be a chain in $\mc_{n-v}$.  
Then $f(x_{p_{h-1}},y_{p_h})=1$ for $h=2,\dots,b$;  
$f(x_{p_h},y_{p_h})=0$ for $h=1,\dots,b$;  
and $f(y_{p_1},t)=f(s,x_{p_b})=0$.
By~(\ref{eqn:defcri2}), 
$f(s,x_{p_b})=0$ and $\pi (s)=0$ imply  
$\pi (x_{p_b})\leq 0$, so $\pi (x_{p_b})=0$. 
Similarly, $f(y_{p_1},t)=0$ and $\pi (t)=p$ imply  
$\pi (y_{p_1})\geq p$ and thus $\pi (y_{p_1})=p$.  
Also, 
$f(x_{p_h},y_{p_h})=0$ and (\ref{eqn:defcri2}) 
imply $\pi (y_{p_h})\leq 1+\pi (x_{p_h})$, for $h=1,\dots,b$. 
Finally, 
$f(x_{p_{h-1}},y_{p_h})=1$ and (\ref{eqn:defcri1}) 
imply $\pi (y_{p_h})\geq \pi (x_{p_{h-1}})$,
for $h=2,\dots,b$. 

The last two statements mean that the sequence 
\[
0=\pi (x_{p_b}),\pi (y_{p_b}),\pi (x_{p_{b-1}}),
\pi (y_{p_{b-1}}),\dots,\pi (y_{p_1})=p
\] 
may only increase in increments of~1, and these may only occur at 
steps of the form $(\pi (x_c),\pi (y_c))$, for $c\in C$. 
Therefore, for any $i\in\{1,\dots,p\}$, 
there exists an element $c\in C$ such that 
$\pi (y_c)=i>\pi (x_c)$; if there are several such $c\in C$, let us
take the greatest one. 
Note that $c\in P_i\,$, and  suppose that $c\notin A_i\,$. 
Let $c'\in P_i$ be such that $c'>c$. 
Then $\pi (y_{c'})=\pi (y_c)=i>\pi (x_c)$,
and (\ref{eqn:defcri2}) implies that $f(x_c,y_{c'})=1$. 
However, the latter means that $c'\in C$, contradicting the choice
of~$c$. 
Hence $c\in A_i\,$, and $C$ intersects all antichains $A_i$
of~$\ma_p\,$. 
\endproof

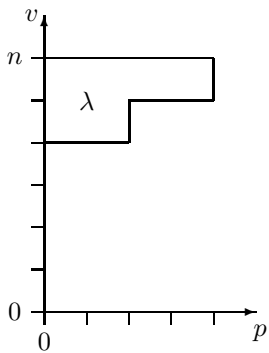
\begin{figure}[ht]
\centering
\setlength{\unitlength}{1.6pt} 
\begin{picture}(60,80)(0,0)
\put(10,10){\vector(1,0){50}}
\put(10,10){\vector(0,1){70}}
\multiput(10,10)(10,0){5}{\line(0,-1){3}}
\multiput(10,10)(0,10){7}{\line(-1,0){3}}
\put(10,3){\makebox(0,0){0}}
\put(3,10){\makebox(0,0){0}}
\put(3,70){\makebox(0,0){$n$}}
\put(61,5){\makebox(0,0){$p$}}
\put(7,80){\makebox(0,0){$v$}}
\put(20,60){\makebox(0,0){$\dc$}}
\put(10,70){\line(1,0){40}}
\thicklines
\put(10,50){\line(1,0){20}}
\put(30,50){\line(0,1){10}}
\put(30,60){\line(1,0){20}}
\put(50,60){\line(0,1){10}}
\put(10,10){\line(0,1){40}}
\end{picture}
\caption{The points $(p,v)$ define the shape $\dc=\dc(P)$}
\label{fig:domain}
\end{figure}

We are now prepared to prove Theorem~\ref{thm:dual}. 
Let us take a closer look at what happens in the course of
Algorithm~\ref{alg:MC}. 
We start with zero potential and zero flow. 
Each step of the algorithm raises the value of the potential or
the flow (but not both). 
The pairs $(p,v)$ occurring in the course of the
algorithm can be represented by points on the coordinate plane;
let us connect these points in the order in which they were
obtained. 
As an example, Figure~\ref{fig:domain} shows the result of applying this
procedure to the partially ordered set $P$ in Figure~\ref{fig:dual}. 
Since both $v$ and $p$ weakly increase during the execution of the
algorithm, the line that connects the points $(p,v)$ defines a Ferrers
shape~$\dc$ (see Figure~\ref{fig:domain}). 
We will prove the Duality Theorem by showing
that the row and column lengths of this shape are exactly the parameters
$\dc_k$ and $\tilde\dc_k$ appearing in Theorem~\ref{thm:dual}. 
(Thus $\dc=\dc(P)$.)
As a byproduct, this will imply that the sequence of points $(p,v)$
generated by the algorithm does not depend on the choice of paths 
used to increase the flow. 

Consider a flow-increasing step 
$(p,v)\leadsto (p,v+1)$
of Algorithm~\ref{alg:MC}, for some $p>0$.
By Theorem~\ref{th:FF}, the flow cost increases by~$p$. 
By the definition of the chain families, 
$|\cup \mc_{n-v}|=|\cup\mc_{n-v-1}|+p$. 
By Lemma~\ref{lem:frank}, 
$\mc_{n-v}$ and $\mc_{n-v-1}$ are orthogonal to 
the antichain family~$\ma_p\,$. 
By Lemma~\ref{lem:A+B=n+kl}, this implies that $\mc_{n-v}$ and $\mc_{n-v-1}$ 
are maximal, and therefore $\dc_{n-v}=p$,
in the notation of Theorem~\ref{thm:dual}. 
Thus each $\dc_k$ is indeed the length of the $k$'th row of~$\dc$
(counting from the top). 

Now consider a potential-increasing step 
$(p,v)\leadsto (p+1,v)$.  
By Lemma~\ref{lem:frank}, the chain family $\mc_{n-v}$ 
is orthogonal to both $\ma_p$ and $\ma_{p+1}$. 
Lemma~\ref{lem:A+B=n+kl} gives
\begin{align*}
\da_{p+1}
&=|\cup \ma_{p+1}|-|\cup \ma_p|
=|\cup \ma_{p+1}|+|\cup\mc_{n-v}|-|\cup \ma_p|-|\cup\mc_{n-v}|\\
&=(n+(p+1)(n-v))-(n+p(n-v))
=n-v \,. 
\end{align*}
Thus $\da_{p+1}$ is the height of the $p+1$'st column of~$\dc$, 
as desired. 
\endproof

The following orthogonality criterion is a direct corollary of
Theorem~\ref{thm:dual} (whose proof has just been completed) and 
Lemma~\ref{lem:A+B=n+kl}. 

\begin{lemma} \label{lem:orth2}
Let $\mc=\{C_1,\dots,C_l\}$ and $\ma=\{A_1,\dots,A_k\}$ 
be chain and antichain families, respectively. 
Then $\mc$ and $\ma$ are orthogonal if and only if the following
conditions hold:
\begin{itemize}
\item[(1)] 
$\mc$ is a maximal chain $l$-family;
\item[(2)]  
$\ma$ is a maximal antichain $k$-family;
\item[(3)] the point $(k,l)$ lies on the outer boundary of the
  shape~$\dc(P)$ (cf.\ Figure~\ref{fig:orth}). 
\end{itemize}
\end{lemma}


\section{Three proofs of Theorem~\ref{thm:mono}} 
\label{sec:proof-mono}

In this section, we provide three proofs of Theorem~\ref{thm:mono},
which use three different ``lines of attack''.
The first proof utilizes Frank's network, and is very much in
the spirit of Section~\ref{sec:frank's-proof}.
The second proof, due to Curtis Greene, 
employs a lattice-theoretic construction
introduced in the original Greene-Kleitman paper~\cite{GK1}.
The third proof, due to Emden Gansner, 
takes advantage of a connection, 
described in Theorem~\ref{thm:Jordan-generic}, 
between posets and linear algebra. 
Although the first proof is somewhat longer and less elegant 
than the second and third proofs, 
it has the advantage of producing, 
as a byproduct, the poset-theoretic analogue of the
augmenting path construction of Ford and Fulkerson
(see Theorem~\ref{thm:fomin1} and Figure~\ref{fig:paths} below). 

The first two proofs will require the following notion. 
Let $\ma$ be an antichain $l$-family in~$P$. 
The \emph{canonical form} of $\ma$ is an antichain 
family $\ma '=\{ A_1 ',\dots,A_l '\}$ 
defined as follows. 
For $i=1,\dots,l$, let $A_i '$ be  
the set of elements $p\in \bigcup \ma$ such that the 
longest chain contained in $\bigcup \ma$ whose top element
is $p$ has length~$i$. 
Each set $A'_i$ is indeed an antichain (possibly empty). 
Since no chain in $\bigcup \ma$ has length $>l$, 
the families $\ma$ and $\ma '=\{ A_1 ',\dots,A_l '\}$ cover 
the same set of elements: $\bigcup \ma '=\bigcup \ma$. 
In particular, if $\ma$ is a maximal antichain $l$-family, 
then so is~$\ma'$. 
We also note that if $p\in A_i '$, $p'\in A_j '$, and $p>p'$,
then~$i>j$.
If $\ma =\ma '$, then $\ma$ is said to be \emph{of canonical form}. 

\subseqnoref{First proof of Theorem~\ref{thm:mono}}
Suppose we are given a chain $k$-family $\mc$ 
and an antichain $l$-family $\ma$ of canonical form, 
and furthermore $\mc$ and $\ma$ are orthogonal. 
Such an orthogonal pair $(\mc,\ma)$ defines 
a flow $f(x,y)$ and a potential $\pi (x)$ 
on Frank's network as follows. 
For each element $p\notin \cup \mc$
(thus $p\in \cup \ma$),  
we set 
$f(s,x_p)=f(x_p,y_p)=f(y_p,t)=1$. 
For each chain $C\in \mc$ and each element $p\in C$ 
which is not minimal in $C$, let $p'\in C$ be the element
covered by~$p$ within~$C$, 
and set 
$f(s,x_p)=f(x_p,y_{p'})=f(y_{p'},t)=1$.
Set the flow along all remaining edges to~0. 
Let us now define the potential. 
Set $\pi (s)=0$ and $\pi (t)=l$. 
For each element $p$ contained in some antichain $A_i\in \ma$,  
set $\pi (x_p)=i-1$ and $\pi (y_p)=i$. 
Now let $p\notin \cup \ma$
(thus $p\in \cup \mc$), and 
suppose that $p\in C\in \mc$. 
Let $p_1<p_2<\cdots<p_l$ be the $l$ elements in $C\cap (\cup \ma)$. 
If $p<p_1\,$, then set $\pi (x_p)=\pi (y_p)=0$. 
If $p>p_l\,$, then set $\pi (x_p)=\pi (y_p)=l$. 
Otherwise, $p_i<p<p_{i+1}$ for some~$i$,
and we set $\pi (x_p)=\pi (y_p)=i$. 
In other words, for any chain $C=(p_1<p_2<\cdots)\in \mc$, 
the sequence of potentials 
$\pi (x_{p_1}),\pi (y_{p_1}),\pi (x_{p_2}),\pi (y_{p_2}),\dots$
has the form $0,\dots,0,1,\dots,1,2,\dots,l-1,l,\dots,l$,  
where all increases are of magnitude~$1$,
and occur between the values 
$\pi (x_p)$ and $\pi (y_p)$ 
with $p\in C\cap (\cup \ma)$. 

We note that the value of the flow $f$ is $n-k$, 
and the value of the potential $\pi$ is~$l$. 

\begin{lemma} \label{lem:defflow} 
For any orthogonal pair $(\mc,\ma)$, with $\ma$ of canonical form, 
the flow $f$ and potential $\pi$ defined as above satisfy the 
conditions {\rm (\ref{eqn:defcri0})-(\ref{eqn:defcri2})}.
\end{lemma}

\proof 
The only nontrivial task is to verify the condition (\ref{eqn:defcri2})
for the edges of the form  $(x_p,y_{p'})$, $p>p'$. 
For such an edge with no flow,
we need to show that $\pi (y_{p'})-\pi (x_p)\leq a(x_p,y_{p'})=0$. 
This is trivially true if 
$\pi (x_p)=l$, or $\pi (y_{p'})=0$, 
or if $p$ and $p'$ are contained in the same chain. 
Assume that none of these statements are true. 
Of the remaining cases, we shall only examine one, 
as the others are dealt with similarly. 
Suppose $p\in C-\cup \ma$ and $p'\in C'-\cup \ma$,
where $C,C'\in\mc$. 
As $\pi (x_p)<l$, 
there is an element $p_0\in C\cap (\cup \ma)$, 
$p_0>p$, 
such that $\pi (x_{p_0})=\pi (x_p)$. 
Similarly, 
there is an element $p_0 '\in C'\cap (\cup \ma)$, 
$p_0 '<p'$, 
such that $\pi (y_{p_0 '})=\pi (y_{p'})$. 
Then $p_0>p>p'>p_0 '$. 
Since $\ma$ is of canonical form, 
we have
$\pi(x_{p_0 '})<\pi (x_{p_0})$ and 
$\pi (y_{p'})-\pi (x_p)
=\pi (y_{p_0 '})-\pi (x_{p_0})
=\pi (x_{p_0 '})-\pi (x_{p_0})+1\leq 0$. 
\endproof

We will now show that in the case of Frank's network,
the minimal cost algorithm may in some sense be reversed.
Suppose we are given a flow $f$ and a potential $\pi$ 
that satisfy the conditions (\ref{eqn:defcri0})-(\ref{eqn:defcri2}). 
The following algorithm iteratively modifies $f$ and $\pi$ 
(hence their respective values $v$ and $p$)
so that at each iteration, either the flow or the potential 
is modified (but not both), and the corresponding value 
(i.e., $v$ or~$p$) decreases by~$1$.  

\begin{algorithm}[Reverse minimal cost flow algorithm for Frank's network]
\label{alg:RMC} {\rm 
\ \par 
\nopagebreak[4]
\medskip
\noindent
\emph{Initial data:} 
flow $f$ and potential $\pi$ in Frank's network satisfying  
{\rm (\ref{eqn:defcri0})--(\ref{eqn:defcri2})}. 

\subseqnoref{RMC1} 
Let $G'=(V,E')$ be the directed graph on the same set of vertices~$V$,
with 
\begin{align}
\notag 
E' =   \ \  & \{(x,y):(x,y)\in E,\ \pi (y)-\pi (x)=a(x,y),\ f(x,y)>0\}\\
\notag \cup &\{(y,x):(x,y)\in E,\ \pi (y)-\pi (x)=a(x,y),\ f(x,y)<c(x,y)\}\,. 
\end{align}
Let $X\subseteq V$ be the set of vertices $x$
for which a path from $s$ to $x$ exists in~$G'$. 
If $t\in X$, then go to {\bf RMC2a}. 
Otherwise go to {\bf RMC2b}.

\subseqnoref{RMC2a} 
Let $\mathcal M$ be a path in $G'$ from $s$ to~$t$. 
Decrease the flow $f$ along $\mathcal M$ by~1. 
Decrease $v$ by~$1$. 
Proceed to {\bf RMC3}.

\subseqnoref{RMC2b} 
Decrease the potential $\pi (x)$ of each vertex $x\in V-X$ by~1. 
Decrease $p$ by~$1$. 
If any potential $\pi (x)$ is equal to $-1$, reset it to~0. 
If any potential $\pi (x)$ is equal to $p+1$, reset it to~$p$. 

\subseqnoref{RMC3} 
If $p>0$ and $v>0$, then  
return to {\bf RMC1}; otherwise stop. 
} 
\end{algorithm} 

The following lemma is a counterpart of Theorem~\ref{th:FF} 
for the reverse minimal cost algorithm.  
Note that in this section,
we only work with Frank's network; 
for general networks, some of the assertions below would be false. 

\begin{lemma} \label{lem:RMC}
At any stage $(p,v)$ of Algorithm~\ref{alg:RMC}, 
the flow $f$ and the potential $\pi$ satisfy the
conditions {\rm (\ref{eqn:defcri0})--(\ref{eqn:defcri2})}. 
\end{lemma}

\proof 
Condition (\ref{eqn:defcri0}) is obviously satisfied
(see {\bf RMC2b}). 
Assume then that $f$ and $\pi$ satisfy conditions 
(\ref{eqn:defcri1})--(\ref{eqn:defcri2}). 
If we now decrease the flow, 
it will be along a path with edges $(x,y)$ 
such that $\pi (y)-\pi (x)=a(x,y)$, 
so (\ref{eqn:defcri1}) and (\ref{eqn:defcri2}) will still hold. 
Suppose then that we cannot find a breakthrough path, 
and thus must decrease $\pi (x)$ by~$1$ for some elements $x \in V$, 
obtaining the new potential function~$\pi'$.
Let us check  (\ref{eqn:defcri1}) for~$f$ and~$\pi'$. 
Suppose $\pi '(y)-\pi '(x)<a(x,y)$;
we need to show that $f(x,y)=0$. 
First we note that 
\begin{equation}
\label{eq:pi(y)-pi(x)}
\pi (y)-\pi (x)\leq \pi '(y)+1-\pi '(x)\leq a(x,y)\,.
\end{equation}
If $\pi (y)-\pi (x)<a(x,y)$, then $f(x,y)=0$, per assumption.
Assume that $\pi (y)-\pi (x)=a(x,y)$. 
By~(\ref{eq:pi(y)-pi(x)}), we then have
$\pi (y)=\pi '(y)+1$ and $\pi (x)=\pi '(x)$. 
Note that $\pi (y)>0$ and $\pi (x)<\pi (t)$. 
If $x=s$, 
then $\pi (y)=\pi (s)+a(s,y)=0$, 
a contradiction. 
If $y=t$, then $\pi (x)=\pi (t)-a(x,t)=\pi (t)$, 
also a contradiction. 
Hence 
$x=x_p$ and $y=y_{p'}$ 
for some elements $p,p'\in P$. 
If $x\notin X$, 
then $\pi (x)=0$. 
Since $\pi (x)-\pi (s)=a(s,x)$ and $(s,x)\notin E'$, 
we have $f(s,x)=0$ and therefore $f(x,y)=0$, by flow conservation. 
If $y\in X$, 
then $\pi (y)=\pi (t)$. 
As $\pi (t)-\pi (y)=a(y,t)$ and $(y,t)\notin E'$, 
we have $f(y,t)=0$ and therefore $f(x,y)=0$. 
If $x\in X$ and $y\notin X$, 
then $(x,y)\notin E'$, which together with
$\pi (y)-\pi (x)=a(x,y)$ implies $f(x,y)=0$, as desired. 
Let us now check  (\ref{eqn:defcri2}) for~$f$ and~$\pi'$. 
Suppose $\pi '(y)-\pi '(x)>a(x,y)$. 
Then 
\begin{equation}
\label{eq:pi(y)-pi(x)2}
\pi (y)-\pi (x)\geq \pi '(y)-(\pi '(x)+1)\geq a(x,y)\,.
\end{equation}
If $\pi (y)-\pi (x)>a(x,y)$, then $f(x,y)=c(x,y)$, per assumption.
Assume that $\pi (y)-\pi (x)=a(x,y)$. 
By~(\ref{eq:pi(y)-pi(x)2}), we then have 
$\pi (x)=\pi '(x)+1$ and $\pi (y)=\pi '(y)$. 
Note that $\pi (y)<\pi (t)$ and $\pi (x)>0$. 
Then $\pi (y)-\pi (x)=a(x,y)\geq 0$ implies 
$\pi (y)>0$ and $\pi (x)<\pi (t)$, 
so $x\notin X$ and $y\in X$. 
This, in turn, implies $f(x,y)=c(x,y)$, as desired. 
\endproof

At each stage $(p,v)$ of Algorithm~\ref{alg:RMC},  
the flow $f$ and potential $\pi$  define 
a chain family $\mc_{n-v}=\mc(f)$ and 
an antichain family $\ma_p=\ma(\pi)$,  
as described in Section~\ref{sec:frank's-proof}. 

\begin{corollary} 
\label{cor:RMC} 
At any stage $(p,v)$ of Algorithm~\ref{alg:RMC},
the families $\mc_{n-v}$ and $\ma_p$ are orthogonal.
\end{corollary} 

\proof
Immediate from Lemma~\ref{lem:RMC} and 
the proof of Lemma~\ref{lem:frank}. 
\endproof

The \emph{comparability graph} ${\rm Comp}(P)$ of a poset~$P$ is the 
undirected graph whose vertices are the elements of $P$
and whose edges connect comparable elements. 
To any chain $k-$family 
$\mc$ in $P$, 
we associate a set of edges $\overline{\mc}$ in ${\rm Comp}(P)$ by 
\[ \overline{\mc}
= \{ (x,y)\,:\,\text{$x$ covers $y$ in some chain $C\in\mc$} \} \,. 
\]
The set $\overline{\mc}$ is a collection of disjoint paths in 
${\rm Comp}(P)$; if $\mc$ does not contain one-element chains,
then it is uniquely recovered from~$\overline{\mc}$. 

We use the notation 
$
\mathcal{X} \triangle \mathcal{Y} = 
(\mathcal{X} \setminus \mathcal{Y}) \cup
(\mathcal{Y} \setminus \mathcal{X})
$
for the symmetric difference of two sets  
$\mathcal{X}$ and $\mathcal{Y}$.  

The following result can be viewed as the poset analogue 
of the Ford-Fulkerson theorem.

\begin{theorem} \label{thm:fomin1} \cite{fomin1}
For any maximal chain $k$-family $\mc$ with $\bigcup\mc\!\neq\!P$
(resp., $k\!>\!0$), 
there exists a maximal chain $(k\!+\!1)$-family~$\mc'$ 
(resp., maximal chain $(k\!-\!1)$-family~$\mc'$) 
such that one of the following is true:
\begin{itemize}
\item
for some path $\mathcal{M}$ in the comparability graph ${\rm Comp}(P)$,  
we have $\overline{\mc'}= \overline{\mc}\triangle \mathcal{M}$,
as shown in Figure~\ref{fig:paths}a
(resp., Figure~\ref{fig:paths}b);
\item
$\mc'$ is obtained from $\mc$ by adding (resp., removing)
a single-element chain. 
\end{itemize}
\end{theorem}

\begin{figure}[ht]
\centering 
\setlength{\unitlength}{1.7pt}
\mbox{
    \subfigure[]{
        \begin{picture}(40,42)(0,0)
        \thicklines
        \put(18,40){$\mc$}
        \put(8,0){\line(0,1){42}}
        \put(20,6){\line(0,1){30}}
        \put(28,0){\line(0,1){42}}
        \put(42,5){\line(0,1){32}}
        \put(-4,16){$\mathcal{M}$}
        \thinlines
        \dashline[+40]{3}(1,24)(9,12)(9,30)(21,12)(21,30)(29,18)(29,30)(37,18) 
        \end{picture}
        } 
\hspace{1in}
    \subfigure[]{
        \begin{picture}(41,42)(0,0)
        \thicklines
        \put(19,40){$\mc$}
        \put(0,24){\line(0,1){18}}
        \put(16,0){\line(0,1){42}}
        \put(28,0){\line(0,1){18}}
        \put(38,5){\line(0,1){28}}
        \put(3,18){$\mathcal{M}$}
        \thinlines
        \dashline[+40]{3}(1,24)(1,36)(17,12)(17,24)(29,6)(29,18)
        \end{picture}
        } 
}
\caption{Augmenting paths in the comparability graph}
\label{fig:paths}
\end{figure}
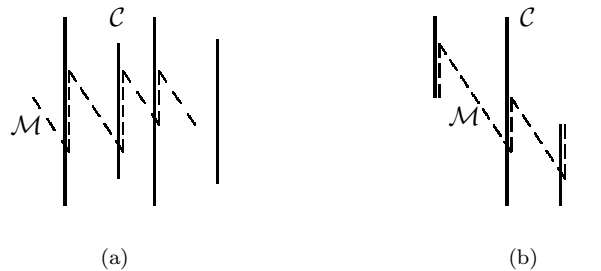 

\proof 
Let $\mc$ be a maximal chain $k-$family. 
By Theorem~\ref{thm:dual}, we can find 
a maximal antichain $l-$family $\ma$ 
of canonical form  
such that $(l,k)$ lies on the boundary of the shape~$\dc(P)$. 
By Lemma~\ref{lem:orth2}, $\mc$ and $\ma$ are orthogonal, 
so we may define a flow $f$ and a potential $\pi$ on the Frank 
network as previously described. 
By Lemma~\ref{lem:defflow}, 
$f$~and $\pi$ satisfy the conditions 
(\ref{eqn:defcri0})--(\ref{eqn:defcri2}), 
so both Algorithm~\ref{alg:RMC} and 
Algorithm~\ref{alg:MC} can be applied. 
In each case, 
we will at some point increase or decrease the flow by~$1$ 
along some path~$\mathcal{N}$.  
Lemma~\ref{lem:frank} (resp., Corollary~\ref{cor:RMC}) implies that 
the resulting families $\mc'$ and $\ma_p$ are orthogonal 
and therefore maximal 
(see Lemma~\ref{lem:orth2}). 
Projecting $\mathcal{N}$ onto~$P$ 
(i.e., applying the map $x_p,y_p\mapsto p$) results in the path $\mathcal{M}$
in $P$ (or a single element $p\in P$)
that can be seen to have the desired properties. 
\endproof

\begin{corollary} 
\label{cor:extremal} 
Let $p$ be an extremal (i.e., minimal or maximal)
element of~$P$. 
Assume that $p\in\bigcup \mc$, 
for any maximal chain $k$-family~$\mc$ (and a fixed~$k$).  
Then $p\in\bigcup \mc'$, for any maximal chain $k'$-family $\mc'$ 
with $k'\geq k$. 
\end{corollary} 

\proof 
It is enough to prove the case $k'=k+1$. 
Let $\mc'$ be a maximal chain $(k+1)$-family. 
By Theorem~\ref{thm:fomin1}, 
there is a path $\mathcal{M}$ in ${\rm Comp}(P)$ 
such that $\overline{\mc}=\overline{\mc'} \triangle \mathcal{M}$ represents 
a maximal chain $k$-family~$\mc$, as shown in Figure~\ref{fig:paths}b
(or else $\mc$ is obtained by removing a single-element chain from~$\mc'$). 
Since $p$ is covered by $\mc$, it must be one of the 
extremal elements of individual chains in~$\mc$. 
It is clear 
that all such elements are also covered by~$\mc'$. 
\endproof

We are now ready to complete the first proof of Theorem~\ref{thm:mono}.  
Let $p$ be an extremal element of $P$. 
Let $k$ be the smallest integer 
such that $p\in\cup \mc$ 
for every maximal chain $k$-family~$\mc$. 
Recall that $c_i(P)$ denotes the number of elements covered by
a maximal chain $i$-family in~$P$. 
Then $c_i(P-\{p\})=c_i(P)$ for all $i=1,\dots,k-1$.
On the other hand, Corollary~\ref{cor:extremal} implies that 
$c_i(P-\{p\})=c_i(P)-1$ for all $i\geq k$.  
Thus the shape $\dc(P-\{p\})$ is identical to~$\dc(P)$, 
except for the $k$th row. 
\endproof

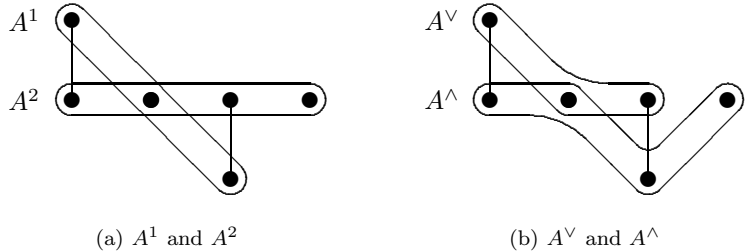
\begin{figure}[ht] 
\centering 
\setlength{\unitlength}{3pt} 
\subfigure[$A^1$ and $A^2$]{
\begin{picture}(36,21)(4,-1)
\put(4,10){\makebox(0,0){$A^2$}}
\put(4,20){\makebox(0,0){$A^1$}}
\multiput(10,20)(10,-10){3}{\circle*{2}} 
\multiput(10,10)(10,  0){4}{\circle*{2}} 
\multiput(10,10)(20,-10){2}{\line(0,1){10}} 
\multiput(10, 8)(0 ,  4){2}{\line(1,0){30}} 
\multiput(8.586,18.586)(2.828,2.828){2}{\line(1,-1){20}}
\put(10,20){\arc(1.414,1.414){180}}
\put(30, 0){\arc(-1.414,-1.414){180}}
\put(10,10){\arc(0, 2){180}} 
\put(40,10){\arc(0,-2){180}} 
\end{picture} 
}
\qquad \qquad 
\subfigure[$A^\vee$ and $A^\wedge$]{
\begin{picture}(36,21)(4,-1)
\put(4,10){\makebox(0,0){$A^\wedge$}}
\put(4,20){\makebox(0,0){$A^\vee$}}
\multiput(10,20)(10,-10){3}{\circle*{2}} 
\multiput(10,10)(10,  0){4}{\circle*{2}} 
\multiput(10,10)(20,-10){2}{\line(0,1){10}} 
\put(10,12){\line(1,0){10}} 
\put(20, 8){\line(1,0){10}} 
\put(8.586,18.586){\line(1,-1){10}} 
\put(21.414,11.414){\line(1,-1){7.172}}
\put(31.414, 4.243){\line(1, 1){7.172}}
\put(31.414,-1.414){\line(1, 1){10}}
\put(10,10){\arc(0, 2){180}}
\put(30,10){\arc(0,-2){180}}
\put(10,20){\arc( 1.414, 1.414){180}}
\put(20,10){\arc(-1.414,-1.414){ 45}}
\put(20,10){\arc( 1.414, 1.414){ 45}}
\put(30, 0){\arc(-1.414,-1.414){ 90}}
\put(30,5.657){\arc(-1.414,-1.414){90}}
\put(40,10){\arc( 1.414,-1.414){180}}
\put(10, 8){\line(1,0){5}} 
\put(25,12){\line(1,0){5}}
\put(11.414,21.414){\line(1,-1){6.4615}}
\put(22.121, 5.047){\line(1,-1){6.4615}}
\put(15,-2.074){\arc(0,10.074){-45}}
\put(25,22.074){\arc(0,-10.074){-45}} 
\end{picture} 
}
\caption{The antichains $A^\vee$ and $A^\wedge$}
\label{fig:GK}
\end{figure}

\subseqnoref{A construction of Greene and Kleitman} 
We will need the following lattice-theoretic construction 
introduced in~\cite{GK1}. 
Let $\ma^1=\{A_1^1,\dots,A_k^1\}$ and $\ma^2=\{A_1^2,\dots,A_k^2\}$ 
be antichain $k$-families of canonical form. 
(Here we allow some of the antichains $A_i^1$ and $A_i^2$ to be empty.) 
Define $\ma^\wedge=\{A_1^\wedge ,\dots,A_k^\wedge \}$ and 
$\ma^\vee =\{A_1^\vee ,\dots,A_k^\vee \}$ by 
\begin{eqnarray}
\label{eqn:Awedge-Avee}
\begin{array}{l}
A_i^\wedge = \{p,\ p\text{ minimal in } A_i^1 \cup A_i^2\}\,,  \\  
A_i^\vee   = \{p \in A_i^1,\exists q\in A_i^2:p \geq q\} \cup
               \{p\in A_i^2,\exists q \in A_i^1:p \geq q\}\,. 
\end{array}
\end{eqnarray} 
To illustrate, consider the antichains $A^1$ and $A^2$ 
in Figure~\ref{fig:GK}a. 
The antichains $A^\wedge$ and $A^\vee$ defined by $A^1$ and $A^2$ 
are shown in Figure~\ref{fig:GK}b. 

\begin{lemma} \label{lem:GK} \cite{GK1} 
$\ma^\vee$ and $\ma^\wedge$ are antichain $k$-families. 
Furthermore, for each $i\leq k$,  we have 
$A_i^\wedge \cup A_i^\vee = A_i^1\cup A_i^2$ and 
$A_i^\wedge \cap A_i^\vee = A_i^1\cap A_i^2$. 
\end{lemma} 

\proof 
For a fixed $i$, a close examination of the definition
(\ref{eqn:Awedge-Avee}) is sufficient to check the second part 
of the lemma, and verify that each $A_i^\wedge$ 
(resp., $A_i^\vee$) is indeed an antichain.
It then remains to prove that the antichains 
$A_1^\wedge,\dots,A_k^\wedge$ (resp., $A_1^\vee,\dots,A_k^\vee$)
are disjoint. 

Suppose $p\in A_i^\wedge \cap A_{j}^\wedge$, 
for $i<j$. 
Assume that $p\in A_{j}^1$. 
Since $\ma^1$ is of canonical form, there exists an element $q\in A_i^1$ 
such that $q<p$. 
This contradicts $p$ being minimal in $A_i^1\cup A_i^2$. 
Similarly, $p\notin A_{j}^2$, 
so $A_i^\wedge$ and $A_{j}^\wedge$ are disjoint. 
Now suppose $p\in A_i^\vee \cap A_{j}^\vee$, $i<{j}$. 
Then either $p\in A_i^1\cap A_{j}^2$ or $p\in A_i^2\cap A_{j}^1$. 
Assume the former. 
By (\ref{eqn:Awedge-Avee}), 
there is an element $q\in A_{j}^1$ such that $p>q$, 
a contradiction since $i<j$ and $\ma^1$ is of canonical form. 
Similarly, $p\notin A_i^2\cap A_{j}^1$,
so $A_i^\vee$ and $A_{j}^\vee$ are disjoint. 
\endproof

\begin{corollary} \label{cor:GK}  \cite{GK1} 
If $\ma^1$ and $\ma^2$ are maximal antichain $k$-families,
then $\ma^\vee$ and $\ma^\wedge$ are maximal as well. 
More generally, if $\ma^1$ is a maximal antichain $k$-family, 
and 
$\ma^2=\{A_1^2,\dots,A_l^2,\underbrace{\phi,\dots,\phi}_{k-l}\}$,
where $\{A_1^2,\dots,A_l^2\}$ is a maximal antichain $l$-family,
$l\leq k$,
then $\ma^\wedge$ is a maximal antichain $k$-family,
and $\ma^\vee=\{A_1^\vee,\dots,A_l^\vee,\underbrace{\phi,\dots,\phi}_{k-l}\}$,
where $\{A_1^2,\dots,A_l^2\}$ is a maximal antichain $l$-family. 
\end{corollary}

\proof
   From Lemma~\ref{lem:GK} 
and the inclusion-exclusion formula, we obtain 
$|\bigcup \ma^\vee |+|\bigcup \ma^\wedge |
=|\bigcup \ma^1|+|\bigcup \ma^2|=a_k(P)+a_l(P)$,
and the claim follows. 
\endproof
 
\subseqnoref{Second proof of Theorem~\ref{thm:mono}}
This proof was contributed by Curtis Greene (private communication). 
The statement below is the antichain analogue of Corollary~\ref{cor:extremal}.

\begin{proposition} 
\label{prop:extremal-antichain} 
Let $p$ be an extremal (i.e., minimal or maximal)
element of~$P$. 
Assume that $p\in\bigcup \ma$, 
for any maximal antichain $k$-family~$\ma$ (and a fixed~$k$).  
Then $p\in\bigcup \ma^1$, for any maximal antichain $k'$-family $\ma^1$ 
with $k'\geq k$. 
\end{proposition} 

\proof
It is enough to consider the case where $p$ is minimal.
Let $\ma^1$ be a maximal antichain $k'$-family of canonical form, 
and suppose that $p\notin \bigcup\ma^1$. 
Let $\ma^2=\{A_1^2,\dots,A_k^2,\underbrace{\phi,\dots,\phi}_{k'-k}\}$,
where $\{A_1^2,\dots,A_k^2\}$ is a maximal antichain $k$-family of canonical form. 
We know that $\ma^2$ covers $p$, and therefore $p\in A_1^2$
because $p$ is minimal and $\ma^2$ is of canonical form. 
Then (\ref{eqn:Awedge-Avee}) and the minimality of $p$ give $p\in A_1^\wedge\,$. 
Since $p\notin A_1^1\,$, Lemma~\ref{lem:GK} implies that $p\notin A_1^\vee$. 
Hence $p$ is not covered by $\ma^\vee$---a contradiction, since by 
Corollary~\ref{cor:GK}, $\ma^\vee$
is, up to a few empty antichains, a maximal antichain $k$-family. 
\endproof

The rest of the proof is straightforward:  
we essentially duplicate the last argument of the first proof,
with chains replaced by antichains.
Let $p$ be an extremal element of $P$, and $k$ the smallest integer 
such that every maximal antichain $k$-family covers~$p$. 
Then $a_i(P-\{p\})=a_i(P)$ for $i=1,\dots,k-1$,
while Proposition~\ref{prop:extremal-antichain} implies that 
$a_i(P-\{p\})=a_i(P)-1$ for all $i\geq k$.  
Thus the shape $\dc(P-\{p\})$ is identical to~$\dc(P)$, 
except for the $k$th column, which is one box shorter. 
\endproof

\subseqnoref{Third proof of Theorem~\ref{thm:mono}}
This proof is due to E.~R.~Gansner~\cite{gansner}.  
It rests on Theorem~\ref{thm:Jordan-generic} 
and the following elementary linear-algebraic lemma
(see \cite{gansner} or \cite{leeuwen2} for a proof). 

\begin{lemma} \label{lem:gansner31} 
Let $V$ be an $n$-dimensional complex vector space,
and $T:V\mapsto V$ a nilpotent linear map with invariants 
$n_1\geq n_2\geq \cdots$. 
Let $W$ be an invariant subspace of~$T$. 
If $T$ is viewed as a nilpotent map $T:V/W\mapsto V/W$ with invariants 
$m_1\geq m_2\geq \cdots$, 
then $n_k\geq m_k$ for all $k\geq 1$. 
\end{lemma}

Assume that $p$ is minimal in~$P$. 
Let $M$ be a generic nilpotent element in~$I(P)$,
viewed as a nilpotent linear map in~$V=\mathbb{C}^n$.
The one-dimensional subspace $W$ spanned by $p$ is $M$-invariant,
and the corresponding map $V/W\mapsto V/W$ has the matrix
obtained from $M$ by striking out the row and the column 
labelled by~$p$. 
Theorem~\ref{thm:Jordan-generic} and Lemma~\ref{lem:gansner31} then imply
that 
$\dc_k(P)=n_k\geq m_k=\dc_k(P-\{p\})$ 
for all $k\geq 1$.
\endproof

\section{Proof of Theorem~\ref{thm:fomin3}
} 

It will be convenient to assume that $p_1,\dots,p_k$ is the  
complete list of \emph{minimal} (rather than maximal) elements of~$P$;
the resulting statement is equivalent to Theorem~\ref{thm:fomin3}
if we pass to the dual poset. 
Assume that $\dc(P-\{p_1\})=\cdots=\dc(P-\{p_k\})=\dc'$. 
The shape $\dc(P)$ is obtained 
by adding a box to~$\dc'$;  
say, this box lies in row $r$ and column~$s$. 
We need to show that $r=k$. 
The number of elements covered by a maximal chain 
$r$-family decreases by~$1$ if any of the $p_i$ is removed from~$P$. 
Hence any maximal chain $r$-family in $P$ covers all the~$p_i\,$,
implying~$r\geq k$. 

Let $\ma=\{ A_1,A_2,\dots\}$ be a maximal antichain $s$-family in~$P$,
and furthermore assume that $\ma$ is of canonical form.
Since the number of elements covered by such a family decreases
if any of the $p_i$ is removed, we conclude that all the $p_i$
are covered by~$\ma$---and therefore contained in~$A_1\,$.
Since any element of $P$ is comparable to some of the~$p_i$, 
the antichain $A_1$ may not contain any other elements,
and its cardinality is equal to~$k$. 
On the other hand, by Lemma~\ref{lem:orth2}, 
$\ma$ is orthogonal to any maximal chain $r$-family,
and therefore any antichain in $\ma$ (including~$A_1$)
must contain at least $r$ elements.
Thus $k\geq r$, and we are done. 
\endproof

\section{Proof of Theorem~\ref{thm:fominG}} 
\label{sec:proof-fominG}

The proof will rely on the Greene-Kleitman construction described
in Section~\ref{sec:proof-mono} (see (\ref{eqn:Awedge-Avee})). 
Suppose the box~$A$ (resp.,~$B$) is located in column~$x_A$ and row~$y_A$
(resp., column~$x_B$ and row~$y_B$). 
In this notation, Theorem~\ref{thm:fominG} 
is equivalent to the following three statements.

\subseqnoref{\ref{thm:fominG}.1} 
If both $p_1$ and $p_2$ are minimal in $P$, then $x_B\leq x_A\,$.

\subseqnoref{\ref{thm:fominG}.2} 
If $p_1$ is maximal and $p_2$ is minimal, 
then $x_B\geq x_A$.

\subseqnoref{\ref{thm:fominG}.3} 
If $p_1$ is maximal and $p_2$ is minimal, 
then $x_B\!\leq\!x_A\!+\!1$.

\subseqnoref{Proof of \ref{thm:fominG}.1} 
Suppose $x_B>x_A\,$. 
On removing $p_2$ from $P$, the value of $a_{x_A}$ does not decrease, 
so there is a maximal antichain $x_A$-family $\ma^1=\{A_1^1,\dots,A_{x_A}^1\}$ 
in $P$ which does not cover~$p_2\,$. 
On~removing $p_1$ from $P-\{p_2\}$, the value of $a_{x_A}$ decreases by~$1$, 
so $p_1$ is covered by any maximal antichain $x_A$-family 
in~$P-\{p_2\}$, including~$\ma^1$. 
Similarly, we may find a maximal antichain $x_A$-family
$\ma^2=\{A_1^2,\dots,A_{x_A}^2\}$ in $P$ that covers $p_2$ but not~$p_1\,$. 
We may assume that $\ma^1$ and $\ma^2$ are of canonical form. 
Let $\ma^\wedge$ and $\ma^\vee$ be defined by (\ref{eqn:Awedge-Avee}). 
By~Corollary~\ref{cor:GK}, 
both $\ma^\vee$ and $\ma^\wedge$ are maximal in~$P$. 
Now $p_1\in A_1^1-\bigcup \ma^2$ and $p_2\in A_1^2-\bigcup \ma^1$, 
so neither $p_1$ nor $p_2$ is contained in $A_1^1\cap A_1^2$. 
Then $p_1$ cannot lie in $A_1^\vee$, 
since this would imply the existence of an element $q\in A_1^2$ 
for which $p_1>q$, 
contradicting the minimality of~$p_1$. 
Similarly, $p_2\notin \bigcup \ma^\vee$. 
Since $p_1,p_2\notin\bigcup\ma^\vee$, and $\ma^\vee$ is maximal, it follows 
that $a_{x_A}(P-\{p_1,p_2\})=a_{x_A}(P)$, 
a contradiction. 
Thus $x_B\leq x_A\,$. 
\endproof

\subseqnoref{Proof of \ref{thm:fominG}.2} 
Suppose~$x_B<x_A\,$. 
The removal of either $p_1$ or $p_2$ from $P$ decreases $a_{x_B}$ by~$1$. 
Hence both $p_1$ and $p_2$ are covered by every maximal antichain 
$x_B$-family in~$P$. 
Choose such a family, of canonical form, and remove $p_1$ from it. 
As~$p_1$ is maximal, the resulting family 
$\ma^1 =\{ A^1_1,\dots,A^1_{x_B}\}$ 
is still of canonical form (and is maximal in $P-\{p_1\}$). 
Removing $p_2$ does not further decrease~$a_{x_B}\,$, 
so there is a maximal antichain $x_B$-family 
$\ma^2=\{ A^2_1,\dots,A^2_{x_B}\}$ in $P-\{p_1\}$ 
for which $\bigcup \ma^2$ does not contain~$p_2\,$. 
We may assume that $\ma^2$ is of canonical form. 
Define $\ma^\vee$ and $\ma^\wedge$ by (\ref{eqn:Awedge-Avee}).
By Corollary~\ref{cor:GK}, 
both $\ma^\wedge$ and $\ma^\vee$ are maximal. 
Suppose that $\ma^\vee$ covers $p_2\,$, i.e., ~$p_2\in A_1^\vee$. 
Since $p_2\notin \bigcup \ma^2\,$, 
there exists an element $p\in A^2_1$ such 
that $p_2>p$, contradicting the minimality of $p_2\,$. 
Hence $p_2\notin \bigcup \ma^\vee$. 
Let $A^1_i\in \ma^1$ be the antichain that used to contain $p_1$ 
before it was deleted from~$P$. 
Then $p_1$ is not comparable to any element $p\in A^1_i$. 
Assume that $p_1$ is comparable to some $p\in A_i^\vee -A^1_i$, 
that is, $p_1>p$. 
By definition of $A_i^\vee$, there exists $q\in A^1_i$ such that~$p>q$. 
But then $p_1>p>q$, a contradiction. 
We conclude that $p_1$ may be added to $A_i^\vee$, 
forming a maximal antichain $x_B$-family in~$P$;
furthermore, this family does not cover $p_2\,$, a contradiction. 
Hence $x_B\geq x_A\,$, as desired.
\endproof

%
%

\subseqnoref{Proof of \ref{thm:fominG}.3} 
The proof is by induction on $n=|P|$.
For $n=1,2$, the claim is easily checked. 
We will rely on Theorem~\ref{thm:mono}, 
Theorem~\ref{thm:fomin3}, and the statements~\ref{thm:fominG}.1 
and~\ref{thm:fominG}.2 above. 
Theorem~\ref{thm:mono} will be used implicitly throughout. 
We may also assume that \ref{thm:fominG}.3 holds for all posets 
of cardinality~$<n$. 

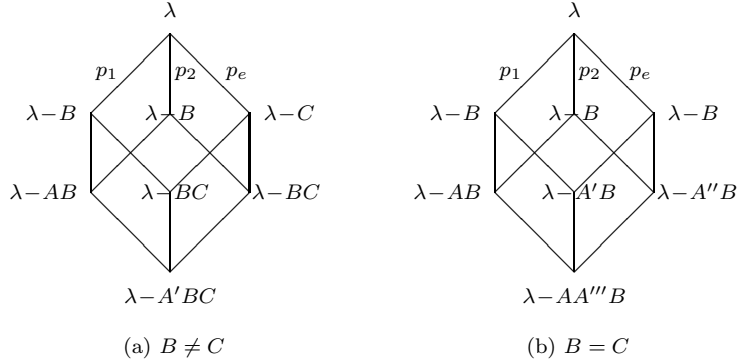
\begin{figure}[ht]
\centering
\setlength{\unitlength}{3pt}
\mbox{
    \subfigure[$B\neq C$]{
    {\footnotesize
        \begin{picture}(40,36)(0,7)
        \multiput(10,20)(20, 0){2}{\line(0, 1){10}}
        \multiput(10,20)( 0,10){2}{\line(1, 1){10}}
        \multiput(10,20)( 0,10){2}{\line(1,-1){10}}
        \multiput(20,10)( 0,10){2}{\line(1, 1){10}}
        \multiput(20,30)( 0,10){2}{\line(1,-1){10}}
        \multiput(20,10)( 0,20){2}{\line(0,1){10}}
        \put(20,43){\makebox(0,0){$\dc$}}
        \multiput(5,30)(15,0){2}{\makebox(0,0){$\dc\!-\!B$}}
        \put(35,30){\makebox(0,0){$\dc\!-\!C$}}
        \put(4,20){\makebox(0,0){$\dc\!-\!AB$}}
        \multiput(20.6,20)(14,0){2}{\makebox(0,0){$\dc\!-\!BC$}}
        \put(20, 7){\makebox(0,0){$\dc\!-\!A'BC$}}  
        \put(12,35){\makebox(0,0){$p_1$}}
        \put(22,35){\makebox(0,0){$p_2$}}
        \put(28.4,35){\makebox(0,0){$p_e$}}  
        \end{picture}
        } }
\qquad
    \subfigure[$B=C$]{
    {\footnotesize
        \begin{picture}(40,36)(0,7)
        \multiput(10,20)(20, 0){2}{\line(0, 1){10}}
        \multiput(10,20)( 0,10){2}{\line(1, 1){10}}
        \multiput(10,20)( 0,10){2}{\line(1,-1){10}}
        \multiput(20,10)( 0,10){2}{\line(1, 1){10}}
        \multiput(20,30)( 0,10){2}{\line(1,-1){10}}
        \multiput(20,10)( 0,20){2}{\line(0,1){10}}
        \put(20,43){\makebox(0,0){$\dc$}}
        \multiput(5,30)(15,0){3}{\makebox(0,0){$\dc\!-\!B$}}
        \put(4,20){\makebox(0,0){$\dc\!-\!AB$}}
        \put(20.6,20){\makebox(0,0){$\dc\!-\!A'B$}}
        \put(35.5,20){\makebox(0,0){$\dc\!-\!A''B$}}
        \put(20, 7){\makebox(0,0){$\dc\!-\!AA'''B$}}
        \put(12,35){\makebox(0,0){$p_1$}}
        \put(22,35){\makebox(0,0){$p_2$}}    
        \put(28.4,35){\makebox(0,0){$p_e$}}
        \end{picture}
        } }
}
\caption{Removing $p_1$, $p_2$, $p_e$ from $P$ 
}
\label{fig:BC}
\end{figure}

We first consider the possibility that $p_1$ (resp.,~$p_2$) is \emph{both} 
minimal and maximal. If that is the case, then any maximal antichain 
contains $p_1$ (resp.,~$p_2$), implying $x_B=1\leq x_A\,$. 
So let us assume that neither $p_1$ nor $p_2$ is both maximal and minimal.
If $p_1$ and $p_2$ are the only extremal elements in~$P$, 
then every element is comparable to both $p_1$ and~$p_2\,$, 
implying $y_A=y_B=1$ and $x_B=x_A+1$.
Let us assume then that $p_e$ is an extremal element different from 
both $p_1$ and~$p_2\,$. 
Denote $\dc-\{C\}=\dc(P-\{p_e\})$. 
The two cases $B\neq C$ and $B=C$ are illustrated in Figure~\ref{fig:BC}.

We first assume that~$B\neq C$.
Let $\dc(P-\{p_1,p_2,p_e\})=\dc-\{A',B,C\}$. 
Since $\dc((P-\{p_e,p_1\})=\dc((P-\{p_e,p_2\})=\dc-\{B,C\}$, 
the induction assumption applied to the poset $P-\{p_e\}$ implies 
$x_B\leq x_{A'}+1$. 
If~$A=A'$, then we are done, so assume $A\neq A'$. 
Then the inclusion $\dc-\{A,B\}\supset \dc-\{A',B,C\}$ implies~$A=C$. 
If~$p_e$ is maximal (and $p_2$ is minimal by assumption), 
then \ref{thm:fominG}.2 applied to the poset~$P-\{p_1\}$ 
implies $x_A\geq x_{A'}$. 
Likewise, if~$p_e$ is minimal, (and $p_1$ is maximal), 
then \ref{thm:fominG}.2 applied to $P-\{p_2\}$ 
implies $x_A\geq x_{A'}$. 
In either case, $x_B\leq x_{A'}+1\leq x_A+1$, as desired. 

It remains to treat the case~$B=C$. 
First assume that $p_e$ is a maximal element. 
As~$p_1$ is also maximal, 
\ref{thm:fominG}.1 applied to $P$ implies $x_B\leq x_{A'}$, 
so if $A=A'$, then $x_B\leq x_A$. 
Assume $A\neq A'$. 
This implies $A'=A'''$. 
If~$A'=A''=A'''$, 
then the induction assumption applied to $P-\{p_e\}$ 
implies $x_{A'}\leq x_A+1$, 
so $x_B\leq x_{A'}\leq x_A+1$. 
The only remaining case, with $B=C$ and $p_e$ maximal, 
is $A'\neq A''=A$. We thus may assume that the latter holds
for any maximal~$p_e\neq p_1\,$. 

The case where $p_e$ is minimal is totally similar. 
As $p_2$ is also minimal, 
\ref{thm:fominG}.1 applied to $P$ implies $x_B\leq x_{A''}$, 
so if $A=A''$, then $x_B\leq x_A$. 
We thus assume~$A\neq A''=A'''$. 
If~$A'=A''=A'''$, 
then the induction assumption applied to $P-\{p_e\}$ 
gives $x_{A'}\leq x_A+1$ and then $x_B\leq x_{A''}\leq x_A+1$. 
The only remaining case, with $B=C$ and $p_e$ minimal, 
is $A''\neq A'=A$. We may furthermore assume that the latter holds
for any minimal~$p_e\neq p_2\,$. 

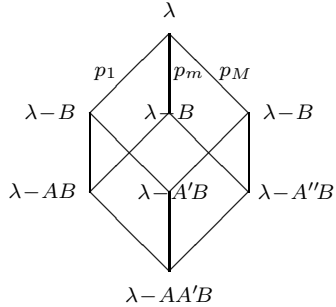
\begin{figure}[ht]
\centering
\setlength{\unitlength}{3pt}
{\footnotesize
\begin{picture}(40,36)(0,7)
\multiput(10,20)(20, 0){2}{\line(0, 1){10}}
\multiput(10,20)( 0,10){2}{\line(1, 1){10}}
\multiput(10,20)( 0,10){2}{\line(1,-1){10}}
\multiput(20,10)( 0,10){2}{\line(1, 1){10}}
\multiput(20,30)( 0,10){2}{\line(1,-1){10}}
\multiput(20,10)( 0,20){2}{\line(0, 1){10}}
\put(20,43){\makebox(0,0){$\dc$}}
\multiput(5,30)(15,0){3}{\makebox(0,0){$\dc\!-\! B$}}
\put(4 ,20){\makebox(0,0){$\dc\! -\! AB$}}
\put(20.5,20){\makebox(0,0){$\dc\! -\! A'\! B$}}
\put(36,20){\makebox(0,0){$\dc\! -\! A''\! B$}}
\put(20, 7){\makebox(0,0){$\dc\! -\! AA'\! B$}}
\put(12,35){\makebox(0,0){$p_1$}}
\put(22.5,35){\makebox(0,0){$p_m$}}
\put(28.4,35){\makebox(0,0){$p_M$}}
\end{picture}
}
\caption{Removing $p_1$, $p_m$, $p_M$ from $P$}
\label{fig:Mm}
\end{figure}

We are now in a situation where for all maximal elements $p_M\neq p_1$ 
and all minimal elements $p_m\neq p_2\,$, we have 
$\dc(P-\{p_m\})=\dc(P-\{p_M\})=\dc-\{B\}$ 
and 
\begin{equation} \label{eqn:ass} 
 \dc(P\!-\!\{p_1,p_M\}) \neq \dc\!-\!\{A,B\} = \dc(P\!-\!\{p_M,p_2\}) 
=\dc(P\!-\!\{p_1,p_m\}). 
\end{equation}
As $\dc(P-\{p_e\})=\dc-\{B\}$ for all extremal elements $p_e$, 
Theorem~\ref{thm:fomin3} implies that $P$ has equally many maximal 
and minimal elements, namely~$y_B\,$. 
If any element $p_e$ is both minimal and maximal, 
then $p_e$ is contained in all maximal antichains, 
so $x_B=1$, which implies $x_B\leq x_A$. 
Assume that no element is both maximal and minimal. 
Then, as there are at least two maximal elements or two minimal elements, 
there exist a maximum element $p_M\neq p_1$ 
and a minimal element $p_m\neq p_2\,$. 
Choose such elements $p_M$ and $p_m\,$. 
By (\ref{eqn:ass}), the situation is as in Figure~\ref{fig:Mm}, 
with $A\neq A'$. 
(The notation $A'$ and $A''$ in Figure~\ref{fig:Mm} is unrelated to 
similar notation in Figure~\ref{fig:BC}.)

Assume~$A\neq A''$. Then $A'=A''$, 
so the induction assumption for~$P-\{p_M\}$ implies $x_{A'}\leq x_A+1$. 
On the other hand,
applying \ref{thm:fominG}.1 to $p_1$ and $p_M$ in~$P$ yields $x_B\leq x_{A'}$, 
so $x_B\leq x_{A'}\leq x_A+1$. 
The only remaining case is $A=A''$.
We may furthermore assume that 
\begin{equation}
\label{eqn:last}
\dc(P-\{p_M,p_m\})=\dc-\{A,B\}
\end{equation}
for \emph{any} maximal element~$p_M$ 
and \emph{any} minimal element $p_m\,$. 
Conditions (\ref{eqn:ass}) show that in (\ref{eqn:last}), the elements 
$p_M$ and $p_m$ do not have to differ from $p_1$ and~$p_2\,$,
respectively. 
Now look at~$P-\{p_2\}$. 
Since $\dc(P-\{p_2,p_M\})=\dc-\{A,B\}$ for all maximal~$p_M$, 
Theorem~\ref{thm:fomin3} implies that $y_A$ is equal to the number of 
maximal elements in $P-\{p_2\}$.
Recall that $y_B$ equals the number of 
maximal elements in~$P$. 
Thus $y_A=y_B$, which implies $x_B=x_A+1$.
\endproof 

\subseqnoref{Sharpness of Theorem~\ref{thm:fominG}} 
We will now demonstrate that Theorem~\ref{thm:fominG} is sharp,
in the sense that its conclusions cannot be strengthened.
Suppose we are given three nested shapes 
$\dc\supset\dc-\{B\}\supset\dc-\{A,B\}$;
thus $B$ is a corner box of $\dc$, while $A$ is a corner box of
$\dc-\{B\}$. 
We need to show that whenever the locations of $A$ and $B$ comply
with one of the two conclusions of  Theorem~\ref{thm:fominG},
there exists a finite poset~$P$ and its extremal
elements $p_1$ and $p_2$ of appropriate type(s) 
such that (\ref{eq:conditionsG}) holds.

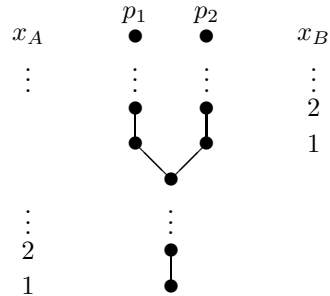
\begin{figure}[ht]
\centering 
\setlength{\unitlength}{2.7pt}
\begin{picture}(40,45)(-5,-5)
\put(-5, 0){\makebox(0,0){$1$}}
\put(-5, 5){\makebox(0,0){$2$}}
\put(-5,10){\makebox(0,0){$\vdots$}}
\put(-5,30){\makebox(0,0){$\vdots$}}
\put(-5,35){\makebox(0,0){$x_A$}}
\put(35,20){\makebox(0,0){$1$}}
\put(35,25){\makebox(0,0){$2$}}
\put(35,30){\makebox(0,0){$\vdots$}}
\put(35,35){\makebox(0,0){$x_B$}}
\put(10,38){\makebox(0,0){$p_1$}}
\put(20,38){\makebox(0,0){$p_2$}}
\put(15,10){\makebox(0,0){$\vdots$}}
\put(10,30){\makebox(0,0){$\vdots$}}
\put(20,30){\makebox(0,0){$\vdots$}}
\put(15, 0){\line(0 ,1){5}}
\put(15,15){\line(-1,1){5}}
\put(15,15){\line( 1,1){5}}
\put(10,20){\line(0 ,1){5}}
\put(20,20){\line(0 ,1){5}}
\multiput(15, 0)(0,5){2}{\circle*{2}}
\put(15,15){\circle*{2}}
\multiput(10,20)(0,5){2}{\circle*{2}}
\multiput(20,20)(0,5){2}{\circle*{2}}
\multiput(10,35)(10,0){2}{\circle*{2}}
\end{picture}
\caption{A subposet used in proving 
sharpness of Theorem~\ref{thm:fominG}}
\label{fig:sub}
\end{figure}

This is easy to do in the cases $x_A=x_B$ 
(regardless of the types of $p_1$ and~$p_2$) and $x_B=x_A+1$, 
as in these cases we may take $P$ to be a disjoint union of chains
whose lengths are the row lengths of~$\dc$. 
The case $x_B<x_A$ is slightly more difficult. 
To construct the poset~$P$, form a disjoint union  
of chains of length $\dc_i\,$, 
for all $i\notin\{y_A,y_B\}$,
together with the $Y$-shaped, $(x_A+x_B)$-element 
subposet in Figure~\ref{fig:sub}. 
Let $p_1,p_2\in P$ be the maximal elements of this subposet. 
The conditions (\ref{eq:conditionsG}) are then easily checked. 
\endproof

\section{Proof of Theorem~\ref{thm:gansner34}}
\label{sec:gansner34}

We use the notation introduced in the first paragraph of
Section~\ref{sec:proof-fominG}. 
Suppose that $x_A\geq x_B\,$. 
Then $y_A<y_B$, implying $c_{y_A}(P-\{p_1\})=c_{y_A}(P)$. 
Let $\mc$ be a maximal chain $y_A$-family in $P-\{p_1\}$. 
Then $p_2$ is contained in some chain $C$ of~$\mc$. 
(Otherwise, removing $p_2$ from $P-\{p_1\}$ would not change the value
of $c_{y_A}\,$.)
Since $C\cup\{p_1\}$ is a chain in $P$,
we obtain a chain $y_A$-family 
$(\bigcup \mc)\cup\{p_1\}$ in $P$ which covers more elements than
$\mc$ does. 
This is a contradiction, since $c_{y_A}(P-\{p_1\})=c_{y_A}(P)$. 
\endproof

\bigskip

\textsc{Acknowledgments.}
This paper would not be written without Gian-Carlo Rota's vigorous
encouragement. We thank Curtis Greene for invaluable advice, 
and for contributing his proof of Theorem~\ref{thm:mono}. 
We also thank Richard Stanley and Andrei Zelevinsky for helpful
comments.

\end{document}